\documentclass[smallextended]{svjour3}                     

\usepackage{color}   
\usepackage{hyperref}
\hypersetup{
    backref=true,
    citecolor=blue,
    colorlinks=true, 
    linktoc=all,     
    linkcolor=blue,  
}

\usepackage{algorithm2e}
\usepackage[latin1]{inputenc}
\usepackage[T1]{fontenc}
\usepackage{bezier}
\usepackage{dsfont}
\usepackage{mathptmx}      
\usepackage{amsfonts}
\usepackage{amsmath}
\usepackage{amssymb}
\usepackage{amsbsy}
\usepackage{amscd}
\usepackage{amsopn}
\usepackage{amssymb}
\usepackage{amstext}
\usepackage{amsxtra}
\usepackage{bm}
\usepackage{graphicx,psfrag}
\usepackage{latexsym}
\usepackage{multirow}
\usepackage{subfig}
\usepackage{etoolbox}

\smartqed  
\DeclareMathAlphabet{\mathcal}{OMS}{cmsy}{m}{n}

\newtoggle{LatexWithProofs}
\togglefalse{LatexWithProofs}

\newtheorem{pLemma}{Lemma}
\newtheorem{pTheorem}{Theorem}
\newcommand{\lpBind}[2]{#1 #2}
\newcommand{\pRef}[1]{(\ref{#1})}
\newcommand{\peLemma}[1]{ \hyperlink{#1}{$\blacktriangle$} \end{pLemma}}
\newcommand{\peTheorem}[1]{ \hyperlink{#1}{$\blacktriangle$} \end{pTheorem}}

\definecolor{darkpink}{rgb}{0.91, 0.33, 0.5}
\definecolor{darksalmon}{rgb}{0.91, 0.59, 0.48}
\definecolor{desertsand}{rgb}{0.93, 0.79, 0.69}
\definecolor{celadon}{rgb}{0.67, 0.88, 0.69}
\definecolor{darkcyan}{rgb}{0.0, 0.55, 0.55}

\iftoggle{LatexWithProofs}
{
\newcommand{\pbClaim}[1]{{\color{red} ?}}
\newcommand{\peClaim}[1]{{\color{red} ?`}}
\newcommand{\pbLemma}[1]{\begin{pLemma} \label{#1}  {\color{darksalmon} #1}}
\newcommand{\pbTheorem}[1]{\begin{pTheorem} \label{#1}  {\color{darksalmon} #1}}

\newcommand{\pbChain}[1]{\[ }
\newcommand{\peChain}[1]{\ {\color{red} \checkmark  \wrm{#1}} \] }

\newcommand{\pePClaim}[1]{ \checkmark}

\newcommand{\pbTClaim}[1]{\begin{equation} \label{#1}  }
\newcommand{\peTClaim}[1]{ \ {\color{red} \checkmark  \wrm{#1}}\end{equation}}
\newcommand{\peEClaim}[1]{ \ {\color{red} \checkmark  \wrm{#1 \bullet} }\end{equation}}

\newcommand{\pbDef}[1]{\begin{equation} \label{#1}  }
\newcommand{\peDef}[1]{ {\color{blue} \bigstar \wrm{#1}} \end{equation}}

\newcommand{\pbHypot}[1]{ \begin{equation} \label{#1}  }
\newcommand{\peHypot}[1]{  \ \ {{\color{blue} \bigstar \wrm{#1}}}  \end{equation} }

\newcommand{\pStop}[1]{{\color{green} \wrm{#1}}}

\newcommand{\pbProof}[2]{\newpage {\color{darksalmon} Proof of {#1} {#2}} \ref{#2} {\hypertarget{#2}{}}}
\newcommand{\peProof}[2]{.\qed{} {\color{darksalmon} End of proof of {#1} {#2}} \ref{#2}}
}
{
\newcommand{\pbChain}[1]{\[}
\newcommand{\peChain}[1]{\]}

\newcommand{\pbClaim}[1]{}
\newcommand{\peClaim}[1]{}
\newcommand{\peTClaim}[1]{ \end{equation}}
\newcommand{\peEClaim}[1]{ \end{equation}}
\newcommand{\pbTClaim}[1]{\begin{equation} \label{#1}}

\newcommand{\peDef}[1]{  \end{equation}}

\newcommand{\lpPlainClaim}[1]{}
\newcommand{\lpEndPlainClaim}[1]{}
\newcommand{\pbDef}[1]{\begin{equation} \label{#1}}

\newcommand{\pbHypot}[1]{\begin{equation} \label{#1}  }
\newcommand{\peHypot}[1]{ \end{equation}}

\newcommand{\pbLemma}[1]{\begin{pLemma} \label{#1}}
\newcommand{\pbTheorem}[1]{\begin{pTheorem} \label{#1}}

\newcommand{\pStop}[1]{\wrm{#1}}

\newcommand{\pbProof}[2]{{Proof of {#1}} \ref{#2}. {\hypertarget{#2}{}}}
\newcommand{\peProof}[2]{.\qed{}}
}
\newcommand{\wsupdx}{\wnorm{\wvec{\hat{x}} - \wvec{x}}_\infty}

\newcommand{\wabs}[1]{\left|#1\right|}

\newcommand{\wcal}[1]{\mathcal{#1}}

\newcommand{\wdx}[1]{{d \! {#1}}}

\newcommand{\wdfdx}[3]{\frac{\partial {#1}}{\partial {#2}}\!\wlr{{#3}}}
\newcommand{\wdfdxx}[3]{\frac{\partial^2 {#1}}{\partial^2 {#2}}\!\wlr{{#3}}}
\newcommand{\wdfdxy}[4]{\frac{\partial^2 \! {#1}}{\partial {#2} \partial {#3}}\!\wlr{{#4}}}

\newcommand{\wdfc}[2]{{#1}'\!\left(#2\right)}
\newcommand{\wdfcc}[2]{{#1}'\!\!\left(#2\right)}
\newcommand{\whessf}[2]{{\nabla^2 \! #1 \left( #2 \right)}}

\newcommand{\wfc}[2]{{#1}\!\left(#2\right)}
\newcommand{\wfcc}[2]{{#1}\!\!\left(#2\right)}

\newcommand{\wfloor}[1]{\lfloor {{#1}} \rfloor }

\newcommand{\wgradf}[2]{{\nabla \! #1}\left(#2\right)}

\newcommand{\wi}[1]{\wrm{i}}

\newcommand{\wlr}[1]{\left( #1 \right)}

\newcommand{\wnorm}[1]{\left\| #1 \right\|}

\newcommand{\wpe}[2]{\wfc{E}{#1\, ;\wvec{\hat{x}}, \wvec{#2}, \wvec{z}}}
\newcommand{\wpeu}[2]{\wfc{E}{#1\, ;\wvec{\hat{x}}, \wvec{#2}, \wvec{z}}}
\newcommand{\wpq}[3]{\wfc{Q}{#1\, ;\wvec{\hat{x}}, \wvec{#2}, \wvec{#3}}}
\newcommand{\wpr}[1]{\wfc{R}{#1\, ;\wvec{\hat{x}}}}
\newcommand{\wpit}[2]{\wfc{P}{#1;\wvec{\hat{x}},\wvec{#2}}}

\newcommand{\wplagr}[2]{ \wfc{\ell_{#1}}{#2;\, \wvec{\hat{x}}}}
\newcommand{\wplagrx}[2]{ \wfc{\ell_{#1}}{#2;\, \wvec{x}}}
\newcommand{\wplagrxu}[3]{ \wfc{\ell_{#1}}{#2;\, \wvec{#3}}}

\newcommand{\wrone}{\mathds R}
\newcommand{\wrn}[1]{{\mathds R}^{#1}}
\newcommand{\wrm}[1]{\mathrm{#1}}

\newcommand{\wrounde}[1]{\wfc{\wrm{fl}}{#1}}

\newcommand{\wset}[1]{{\left\{ #1 \right\}}}

\newcommand{\wst}[1]{\langle{#1}\rangle}

\newcommand{\wvec}[1]{\mathbf{#1}}

\newcommand{\wvone}[1]{\mathds{1}_{#1} }

\begin{document}
\title{The effects of rounding errors in the nodes on barycentric interpolation}
\author{Walter F. Mascarenhas and Andr\'{e} Pierro de Camargo}

\institute{Walter F. Mascarenhas,
supported by grant
2013/10916\-2 from Funda\c{c}\~{a}o
de Amparo \`{a}  Pesquisa do Estado de S\~{a}o Paulo (FAPESP)
, and Andr\'{e} Pierro de Camargo,  supported by grant 14225012012-0 from CNPq.
 \at  Instituto de Matem\'{a}tica e Estat\'{i}stica, Universidade de S\~{a}o Paulo, \\
           Cidade Universit\'{a}ria, Rua do Mat\~{a}o 1010, S\~{a}o Paulo SP, Brazil. CEP 05508-090.\\
              Tel.: +55-11-3091 5411, Fax: +55-11-3091 6134,   \email{walter.mascarenhas@gmail.com}           
}
\maketitle

\begin{abstract}
We analyze the effects of rounding errors in the nodes on polynomial barycentric interpolation. These errors are particularly
relevant for the first barycentric formula with the Chebyshev points of the second kind. Here,
we propose a method for reducing them.
\end{abstract}

\subclass{65D05, 65G50}

\iftoggle{LatexWithProofs}
{
\vspace{1cm}
\hyperlink{at_work}{Latex Proof Version. Not for reading.}\\
}

\section{Introduction}
Given nodes $x_0 < x_1 < \dots   < x_n$, weights $w_0,\dots,w_n$,
and an interval $[x^-,x^+]$,
two formulae for barycentric interpolation of a function
$f: [x^-,x^+] \rightarrow \wrone{}$ are considered.
 The first one is
\begin{equation}
\label{first_formula}
\wfc{p}{x; \, \wvec{x}, \wvec{y},\wvec{w}} := \wlr{ \prod_{k=0}^n \wlr{x - x_k} } \sum_{k=0}^n \frac{w_k y_k}{x - x_k},
\end{equation}
where $y_k = \wfc{f}{x_k}$. The second one is
\begin{equation}
\label{second_formula}
\wfc{q}{x;\wvec{x}, \wvec{y}, \wvec{w}} \ \ := \ \ \left.
\sum\limits_{k = 0}^{n} \frac{w_k y_k}{x - x_k} \right/ \sum\limits_{k = 0}^{n} \frac{w_k}{x - x_k}.
\end{equation}
The first formula is a polynomial for any weights, but it only interpolates $f$ if
\begin{equation}
\label{def_lambda}
w_k = \wfcc{\lambda_k}{\wvec{x}} := \frac{1}{\prod_{j \neq k} \wlr{x_k - x_j}}.
\end{equation}
The second formula always interpolates $f$ at the nodes $x_k$, but it is a polynomial for all $\wvec{y}$
only if the weights are of the form $w_k = \kappa_n \wfc{\lambda_k}{\wvec{x}}$, where $\kappa_n$ is
independent of $k$.

 Historically, Taylor \cite{TAYLOR} and Salzer \cite{SALZER}
 considered weights for which the second formula is a polynomial.
The recent literature is also concerned with strictly rational second formulae
\cite{BERRUT_POLE,BOS,BOSGEN,FLOATER}, and we address this case in \cite{MascCam}.
Here we focus on polynomial interpolation, specially in the
classical case considered by Salzer, in which we interpolate
at the Chebyshev points of the second kind.

In theory, barycentric interpolation at the Chebyshev points leads to accurate results.
In practice, the first barycentric formula with these nodes
suffers from accuracy problems when the number of nodes is large. In order
to understand these problems, we must consider the steps outlined in Figure \ref{figure_steps}:
\begin{center}
\begin{figure}[h]
\begin{picture}(0,85)(0,52)
\put(4,140){Abstract function}
\put(19,129){$\wfc{f}{x}$}
\put(79,140){Step I: Abstract interpolation}
\put(73,137){\vector(1,0){112}}
\put(79,128){Error I: Approximation theory}
\put(193,138){Abstract interpolant}
\put(206,130){$\wfc{a}{x,\wvec{x},\wfc{f}{\wvec{x}},\wvec{w}}$}
\put(269,134){
$
\left\{
\begin{array}{l}
\wvec{x} = \wrm{exact \ nodes}    \\
\wvec{w} = \wrm{exact \ weights}  \\
\end{array}
\right.
$}
\put(227,123){\vector(0,-1){47}}
\put(231,112){Step II: Finite precision}
\put(257,104){representation of $a$.}
\put(231,93){Error II: Given by how $\wvec{x}$, $\wfc{f}{\wvec{x}}$ }
\put(260,85){and $\wvec{w}$ are rounded.}
%
%
\put(207,65){In practice we}
\put(207,57){use $\wfc{a}{x,\wvec{\hat{x}},\wvec{y},\wvec{\hat{w}}}$}
\put(256,61){
$
\left\{
\begin{array}{l}
\wvec{\hat{x}} = \wrm{rounded \ nodes} \\
\wvec{y} = \wrm{rounded} \ \wfc{f}{\wvec{x}} \\
\wvec{\hat{w}} = \wrm{rounded \ weights}
\end{array}
\right.
$}
\put(26,120){\vector(0,-1){42}}
\put(30,108){The overall error is a}
\put(30,100){combination of the}
\put(30,92){errors in the three steps.}
\put(198,64){\vector(-1,0){119}}
\put(88,68){Step III: Evaluation of $\wfc{a}{x,\wvec{\hat{x}},\wvec{y}, \wvec{\hat{w}}}$}
\put(87,55){Error III: Usual stability analysis}
\put(9,66){Final result}
\put(0,56){$\wfc{f}{x} \ \ \approx \ \ \wrounde{\wfc{a}{x,\wvec{\hat{x}},\wvec{y}, \wvec{\hat{w}}}}$}
\put(17,49){$???$}
\end{picture}
\caption{The overall error in interpolation. In this article $a$ is either $p$ in \pRef{first_formula} or $q$  in \pRef{second_formula}.}
\label{figure_steps}
\end{figure}
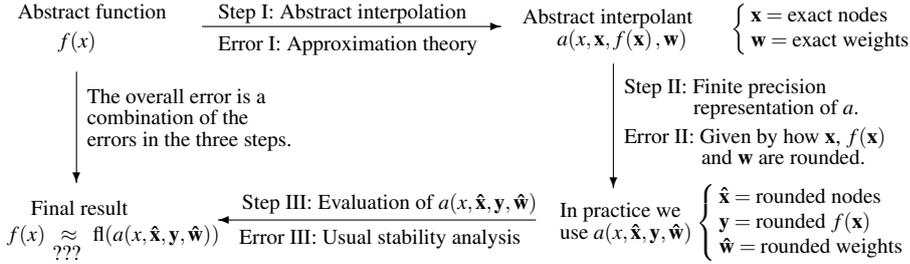
\end{center}

The errors in Step II are not considered
to their full extent in the current literature. For instance,
\cite{HIGHAM_STAB} takes into account the rounding errors in the  evaluation of
\pRef{def_lambda}, but it does not consider weights obtained from closed form
expressions evaluated at exact nodes, as in Salzer's case.
Table \ref{table_rho} compares the errors introduced by Step II and by Step III
for usual implementations of formulae \pRef{first_formula} and \pRef{second_formula} in
Salzer's case over two
sets of points $X_{-1,n}$ and $X_{0,n}$ described in Appendix \ref{section_experiments}.
It shows that the errors introduced by Step II can be larger than those introduced by Step III for both formulae.
The entries in bold face in Table \ref{table_rho}
highlight Step II errors that are much larger than Step III errors for the first formula.

\begin{table}[!h]
\caption{The ratio $\frac{\wrm{mean \ of \ errors \ in \ Step \ II}}{\wrm{mean \ of \ errors \ in \ Step \ III}}$
in  Salzer's case
for the sets $X_{-1,n}$ and $X_{0,n}$.}
\centering
\begin{tabular}{c|cc|cc|cc|cc}
\hline\\[-0.27cm]
      & \multicolumn{4}{c|}{First Formula}  &  \multicolumn{4}{c}{Second Formula} \\[-0.01cm]
\cline{2-9}\\[-0.22cm]
      & \multicolumn{2}{c|}{$\wfc{f}{x} = \cos x$}  &  \multicolumn{2}{c|}{$\wfc{f}{x} = \wfc{\cos}{10^4 x}$}
      & \multicolumn{2}{c|}{$\wfc{f}{x} = \cos x$}  &  \multicolumn{2}{c}{$\wfc{f}{x} = \wfc{\cos}{10^4 x}$}  \\
\cline{2-9}\\[-0.25cm]
$n + 1$    &$x \in X_{-1,n}$ & $x \in   X_{0,n}$ & $x \in  X_{-1,n}$ & $x \in X_{0,n}$ & $x \in  X_{-1,n}$    & $ x \in X_{0,n}$     & $x \in X_{-1,n}$ & $x \in X_0$  \\[-0.01cm]
\hline\\[-0.2cm]
$ 10^3$& {\bf 1.4e+2} & 7.8e-1 & Step I is         & Step I is      & 4.4e-2 & 3.7e-2 & Step I is                 & Step I is \\
$ 10^4$& {\bf 4.4e+3} & 1.0e-1 & critical          & critical       & 1.7e-2 & 1.0e-2 & critical                  & critical\\[0.1cm]
$ 10^5$& {\bf 1.5e+5} & 5.2e-1 & {\bf 1.5e+5} & 2.8e-1 & 5.6e-3 & 3.2e-3 & 6.6 & 8.9e-2\\
$ 10^6$& {\bf 8.4e+6} & 7.1e-2 & {\bf 7.3e+6} & 6.6e-2 & 1.9e-3 & 1.0e-3 & 4.0 & 1.2e-2
\end{tabular}
\label{table_rho}
\end{table}

This article estimates the errors in Steps II and III for the two barycentric formulae and  proposes
a strategy to reduce these errors in practice.
In the next section we
describe an experiment with the sine function which corroborates the
data in Table \ref{table_rho} and present an overview of our results.
In Section \ref{section_notation} we present our notation and estimates for the
order of magnitude of the parameters relevant to our analysis of Salzer's case.
Section \ref{section_polynomials} analyzes how rounding errors in the nodes affect polynomial interpolation.
Section \ref{section_first} estimates the backward and forward errors for the first
formula, and Section \ref{section_second} presents bounds on the forward
errors for the second formula (we do not present bounds on the
backward error for the second formula because it
is discussed in detail in \cite{MascCam}.)
Our strategy for reducing the errors in Step II is presented in the last section.
Appendix \ref{section_proofs}  proves the lemmas and theorems stated in the previous sections
whereas Appendix \ref{section_experiments} describes the numerical experiments on which our tables and figures are based.

\hypertarget{link_here}
Finally, a question of (modern) notation. The statements of our lemmas and theorems
end with a \hyperlink{link_here}{blue $\blacktriangle$}. By clicking on this triangle you will move to
the proof of the corresponding result. Once you are at the proof, you can click
on the lemma or theorem number in order to return to its statement. You can also
use the navigation buttons in your reader in order to return to a previous page,
or use short cuts in some readers. For instance,
in the Adobe Acrobat reader you can type Alt+left arrow in order to return
to the previous view, and you can also install a button in the toolbar for
this purpose.

\section{Overview and motivation}
\label{section_overview}
The complete analysis of the stability of the barycentric formulae
requires much attention to detail, and people guided by concrete examples will have
a better chance of understanding the subtle points.
For this reason, throughout the
article we illustrate the use of our general results in the following
specific situation:\\

{\bf \hypertarget{link_salzer}{Salzer's Case.}}
{\it
We consider floating point nodes $\wvec{\hat{x}}^c$
obtained by rounding the abstract Chebyshev points of  the second kind
$\wvec{x}^c$:
\begin{equation}
\label{cheby_points}
x^{(c)}_k := - \wfc{\cos}{k \pi / n }
\hspace{0.7cm} \wrm{and} \hspace{0.7cm}
\hat{x}^{(c)}_k := \wfc{\wrm{rounded}}{x^{(c)}_k} := \wrounde{x^{(c)}_k}.
\end{equation}
The weights used in computation are given in closed form by \cite{SALZER}.
These weights are equivalent to $\hat{\wvec{w}} = \wfc{\lambda}{\wvec{x}^c}$,
for $\lambda$ in \pRef{def_lambda}, and we call them Salzer's weights.
We make conservative assumptions about $n$ and the magnitude of the rounding errors.
Formally, we suppose that the nodes are rounded as usual and
\pbHypot{salzersCase}
10 \leq n \leq 2 \times 10^6
\hspace{1cm} \pStop{and} \hspace{1cm}
\wsupdx \leq 4.6 \times 10^{-16}.
\peHypot{salzersCase}
}

Salzer's case is relevant first because of its practical importance, second
because it shows clearly that the concern with the perturbation in
the nodes is not futile. In fact, if we neglect
these errors then we can underestimate the errors
in the first formula by orders of magnitude in this case.
Therefore, in order to fully understand the accuracy of the barycentric formulae $p$ and $q$
in \pRef{first_formula} and \pRef{second_formula}, one must
be aware of the differences between their variations (a), (b) and (c) below,
in which $\hat{\wvec{w}}$ are the weights used in computation
(in Salzer's case $\hat{\wvec{w}} = \wfc{\lambda}{\wvec{x}^c}$ and \cite{HIGHAM_STAB} considers
$\hat{\wvec{w}} = \wrounde{\wfc{\lambda}{\hat{\wvec{x}}}}$.)
\begin{equation}
\label{eq_variations}
(\wrm{a}) \
\left\{
\begin{array}{c}
 \wfc{p}{x; \, \wvec{x}, \wvec{y},\wfc{\lambda}{\wvec{x}}} \\[0.07cm]
  =  \\
 \wfc{q}{x; \, \wvec{x}, \wvec{y},\wfc{\lambda}{\wvec{x}}} \\
 = \\
 \wfc{P}{x;\wvec{x},\wvec{y}}
\end{array}
\right.
\hspace{0.3cm}
(\wrm{b})  \
\left\{
\begin{array}{c}
 \wfc{p}{x; \, \wvec{\hat{x}}, \wvec{y}, \hat{\wvec{w}}} \\[0.07cm]
   \neq  \\
 \wfc{q}{x; \, \wvec{\hat{x}}, \wvec{y},\hat{\wvec{w}}} \\
  \neq \\
 \wfc{P}{x;\wvec{x},\wvec{y}} \neq \wfc{P}{x;\wvec{\hat{x}},\wvec{y}}
\end{array}
\right.
\hspace{0.3cm}
(\wrm{c}) \
\left\{
\begin{array}{c}
 \wfc{p}{x; \, \wvec{\hat{x}}, \wvec{y},\wfc{\lambda}{\wvec{\hat{x}}}}\\[0.07cm]
 = \\
  \wfc{q}{x; \, \wvec{\hat{x}}, \wvec{y},\wfc{\lambda}{\wvec{\hat{x}}}} \\
    = \\
 \wfc{P}{x;\wvec{\hat{x}},\wvec{y}},
 \end{array}
\right.
\end{equation}
where
\[
\wfc{P}{x;\wvec{x},\wvec{y}} := \wrm{the \ } n \wrm{th \ degree \ polynomial \ that \ interpolates \ } \wvec{y} \ \wrm{at} \ \wvec{x}.
\]

In variation $(\wrm{a})$ we consider the theoretical nodes,
like the Chebyshev points of the four kinds.
Usually, these nodes cannot be represented exactly in finite
precision arithmetic and
in practice we use rounded nodes instead, and these rounded nodes are
considered in variations $(\wrm{b})$ and $(\wrm{c})$.
The rounded nodes do not have  all the theoretical properties of the exact ones,
like the orthogonality of the corresponding polynomials with respect to
convenient inner products or neat closed form expressions
relating them. Unfortunately,
the advantages given by these theoretical properties may be illusory for large $n$,
and we may be subject to subtle side effects when we apply results derived
for the exact nodes to the rounded ones. For instance, if
we use Fourier transforms based on the exact nodes
to obtain weights for the barycentric formulae then we obtain weights
like Salzer's. However, in the following experiment with the first formula
applied to the sine function, the maximum forward error
\begin{equation}
\label{trial_points}
\max_{\wrm{trial \ points \ } x} \wabs{\ \wfc{\sin}{x} - \wfc{p}{x;\wvec{\hat{x}}^c, \wfc{\sin}{\wvec{\hat{x}}^c}, \wfc{\lambda}{\wvec{x}^c}} \ }
\end{equation}
corresponding to the Salzer's  weights $\wfc{\lambda}{\wvec{x}^c}$
was about 650 times larger than the error corresponding to the numerical
weights $\wrounde{\wfc{\lambda}{\wvec{\hat{x}}^c}}$ for $n = 1000$.

\begin{figure}[!h]
\includegraphics[viewport= -80 0 670 590, width=6.5cm, height=5.0cm]{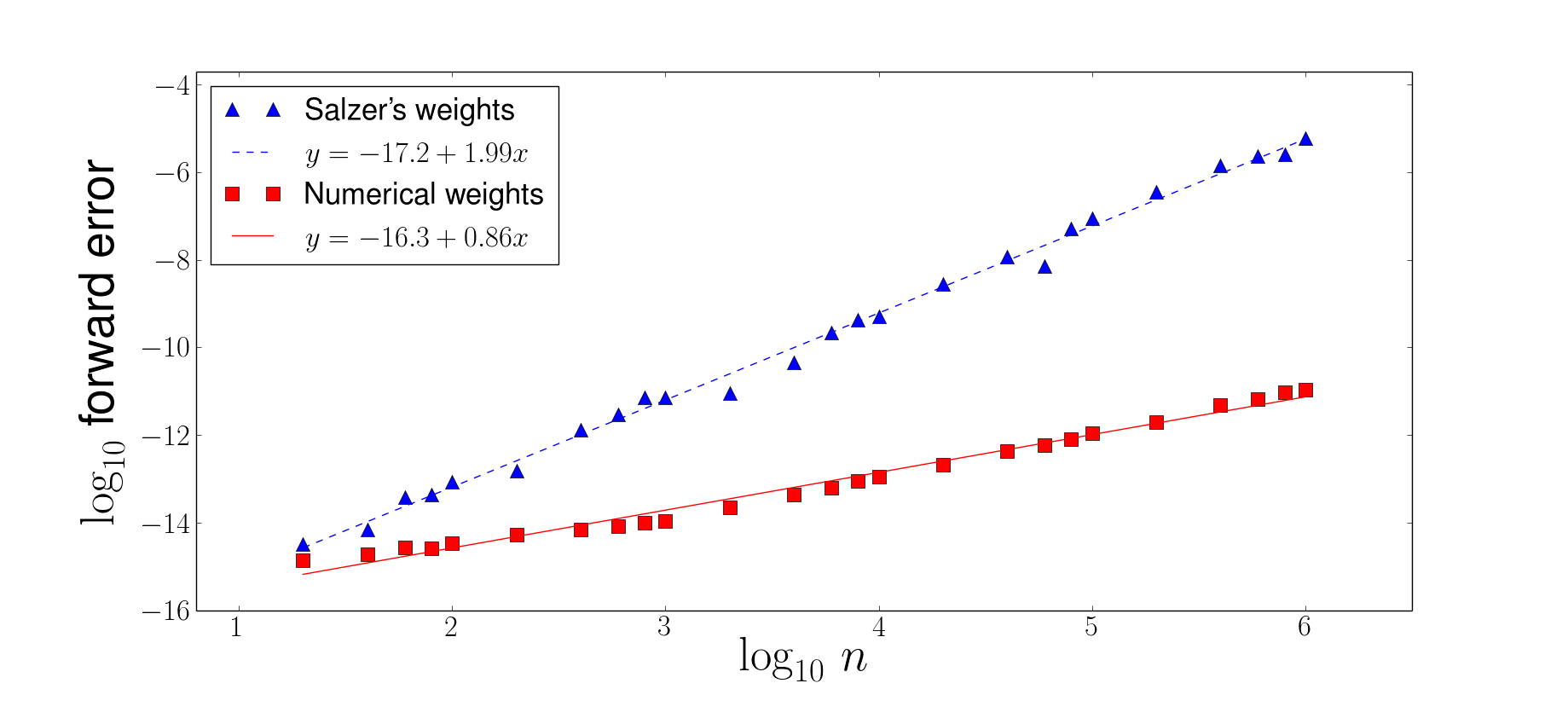}
\caption{$\log_{10}$ of the maximum error in the approximation of $\wfc{f}{x} = \sin x$ in $[-1,1]$ by the first
barycentric formula with rounded nodes $\wvec{\hat{x}}^c$ and (i) Salzer's weights,
$\wfc{\lambda}{\wvec{x}^c}$,
and (ii) the weights corresponding to the rounded nodes, $\wrounde{\wfc{\lambda}{\wvec{\hat{x}}^c}}$, which
we call by {\it Numerical weights}.}
\label{figure_least_squares}
\end{figure}

The set of trial points in \pRef{trial_points}, and
the details of the experiment in Figure \ref{figure_least_squares}, are
presented in appendix \ref{section_experiments}, the nodes $\wvec{\hat{x}}^c$
in this experiment were computed with machine precision $\epsilon \approx 2.3 \times 10^{-16}$
and the straight lines in this plot were obtained by the least squares method.
The straight line for Salzer's weights shows that, in this particular experiment,
the corresponding errors grow like $0.1 \epsilon n^2$ whereas
the errors incurred when using rounded nodes in combination
with the Numerical weights are in better agreement with the
$\wfc{O}{\epsilon n}$ upper bounds presented in \cite{HIGHAM_STAB}.
In Section \ref{section_first} we explain this $\wfc{O}{\epsilon n^2} \times \wfc{O}{\epsilon n}$
discrepancy in the order of magnitude of the forward errors
for the first formula.
In summary, we show that, for large $n$, the maximum forward error for the first formula in Steps II and III in Salzer's case is well
described by $\max \wabs{y_k z_k^s}$, where
$z_k^s := \wlr{\wfcc{\lambda_k}{\wvec{x}^c}- \wfcc{\lambda_k}{\wvec{\hat{x}}^c}}/ \wfcc{\lambda_k}{\wvec{\hat{x}}^c}$.
The experimental evidence shows clearly that $\wnorm{\vec{z}^s}_\infty$ is of order $\epsilon n^2$
in Salzer's case, whereas the analogous $\wnorm{\vec{z}^r}_\infty$ for the rounded nodes in
combination with the Numerical weights is of order $\epsilon n$
(see Tables 1 and 2 in \cite{MascCam}, in which $\zeta_k = - z_k /(1 + z_k) \approx -z_k$.)

For the second formula, we show that the forward error in Step II can be estimated
via the {\it Error Polynomial} $\wpe{x}{y}$, which is given by
\begin{equation}
\label{def_error_pol}
\wpe{x}{y} :=  \wpit{x}{y \wvec{z}} - \wpit{x}{y} \wpit{x}{\wvec{z}},
\end{equation}
where $\wvec{y} \wvec{z}$ is the vector with entries $\wlr{y z}_k = y_k z_k$.
The $z_k$ and the Error Polynomial are the key factors for
the understanding of the maximum forward errors in Step II for the first and
the second barycentric formula in Salzer's case. For this reason, in the following sections we study them
in detail, and explain how they can be bounded in terms of the machine precision $\epsilon$,
\[
L := \wrm{Lipschitz \ constant \ of \ the \ function \ } f \ \wrm{we \ are \ interpolating},
\]
the {\it Lebesgue constant}
\[
\Lambda := \Lambda_{x^-,x^+,\wvec{x}} := \sup_{x \in [x^-,x^+], \ \ \wvec{y} \neq 0} \frac{\wabs{\wfc{P}{x;\wvec{x},\wvec{y}}}}{\wnorm{\wvec{y}}_\infty},
\]
and terms related to the node spacing.
Our conclusions are summarized by the diagram in Figure \ref{figure_bounds}, in which $P^*$ is the best
polynomial approximation of $f$.
\begin{center}
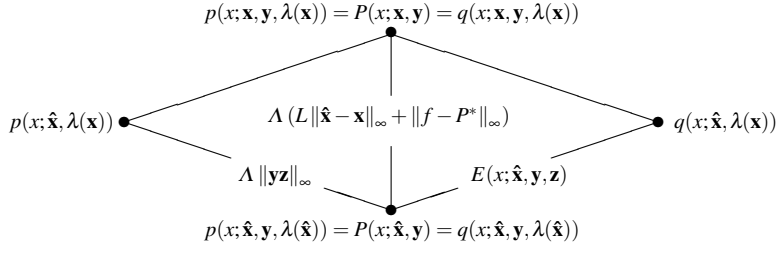
\begin{figure}[h]
\begin{picture}(0,100)(30,5)
\put(100,60){\circle*{4}}
\put(57,58){$\wfc{p}{x;\wvec{\hat{x}},\wfcc{\lambda}{\wvec{x}}}$}
\put(100,60){\line(3,1){100}}
\put(100,60){\line(3,-1){40}}
\put(200,27){\line(-3,1){25}}
\put(143,38){$\Lambda \wnorm{\wvec{y} \wvec{z}}_\infty$}
\put(200,94){\line(0,-1){24}}
\put(154,60){$\Lambda \wlr{L \wnorm{\wvec{\hat{x}} -\wvec{x}}_\infty + \wnorm{f - P^*}_\infty}$}
\put(200,27){\line(0,1){24}}
\put(200,27){\circle*{4}}
\put(200,94){\circle*{4}}
\put(130,99){$\wfc{p}{x;\wvec{x},\wvec{y},\wfcc{\lambda}{\wvec{x}}} =
\wfc{P}{x;\wvec{x},\wvec{y}} = \wfc{q}{x;\wvec{x},\wvec{y},\wfcc{\lambda}{\wvec{x}}}$}
\put(300,60){\circle*{4}}
\put(300,60){\line(-3,1){100}}
\put(300,60){\line(-3,-1){40}}
\put(200,27){\line(3,1){25}}
\put(230,38){$\wfc{E}{x;\wvec{\hat{x}},\wvec{y},\wvec{z}}$}
\put(306,58){$\wfc{q}{x;\wvec{\hat{x}},\wfcc{\lambda}{\wvec{x}}}$}
\put(130,17){$\wfc{p}{x;\wvec{\hat{x}},\wvec{y},\wfcc{\lambda}{\wvec{\hat{x}}}} = \wfc{P}{x;\wvec{\hat{x}},\wvec{y}} =
\wfc{q}{x;\wvec{\hat{x}},\wvec{y},\wfcc{\lambda}{\wvec{\hat{x}}}}$}
\end{picture}
\caption{The forward errors for the formulae $p$ and $q$ in Step II. The values in the edges
of this diagram are upper bounds on the order of magnitude of the maximum difference between their
end points. These bounds hold under technical conditions described in the next sections.}
\label{figure_bounds}
\end{figure}
\end{center}

In the following sections we look carefully at the differences
\[
\wfc{p}{x,\wvec{\hat{x}},\wvec{y},\wfc{\lambda}{\wvec{x}}} - \wfc{p}{x;\wvec{\hat{x}},\wvec{y},\wfc{\lambda}{\wvec{\hat{x}}}}
\hspace{0.6cm} \wrm{and} \hspace{0.6cm}
\wfc{q}{x,\wvec{\hat{x}},\wvec{y},\wfc{\lambda}{\wvec{x}}} - \wfc{q}{x,\wvec{\hat{x}},\wvec{y},\wfc{\lambda}{\wvec{\hat{x}}}}
\]
shown in the lower edges in Figure \ref{figure_bounds}, and justify the estimates presented in these edges.
Using the vertical edge in this figure  and the triangle inequality, we can
bound
\[
\wfc{p}{x,\wvec{\hat{x}},\wvec{y},\wfc{\lambda}{\wvec{x}}} - \wfc{p}{x,\wvec{x},\wvec{y},\wfc{\lambda}{\wvec{x}}}
\hspace{0.6cm} \wrm{and} \hspace{0.6cm}
\wfc{q}{x,\wvec{\hat{x}},\wvec{y},\wfc{\lambda}{\wvec{x}}} - \wfc{q}{x,\wvec{x},\wvec{y},\wfc{\lambda}{\wvec{x}}},
\]
and in Section \ref{section_polynomials} we present bounds on the difference
 $\wfc{P}{x;\wvec{x},\wvec{y}} - \wfc{P}{x;\wvec{\hat{x}},\wvec{y}}$ corresponding to this vertical edge.

\begin{figure}[!h]
\subfloat[$\wpe{x}{y}$]{
\includegraphics[viewport= 60 -20 550 400, width=3.65cm, height=2.0cm]{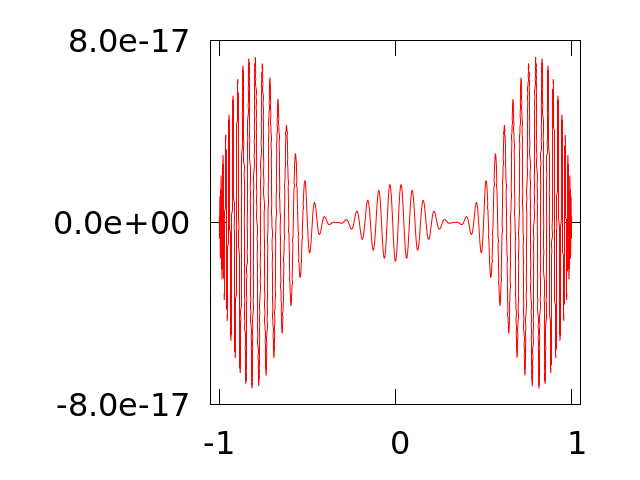}
}
\subfloat[$\wpr{x}$]{
\includegraphics[viewport= -100 10 550 400, width=3.65cm, height=2.0cm]{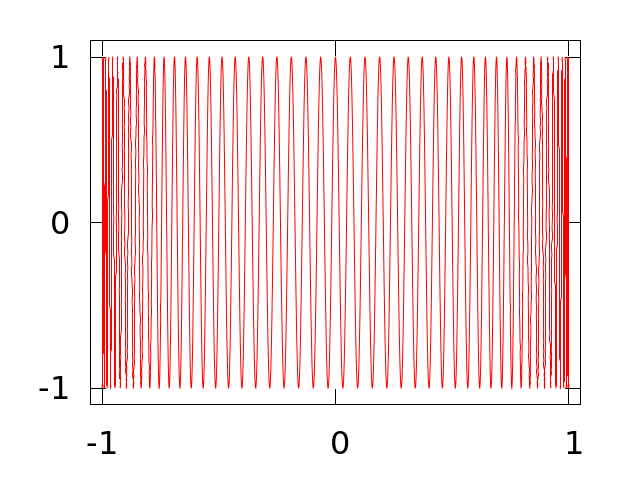}
}
\subfloat[$\wpq{x}{y}{z}$]{
\includegraphics[viewport= -80 10 500 410, width=3.65cm, height=2.0cm]{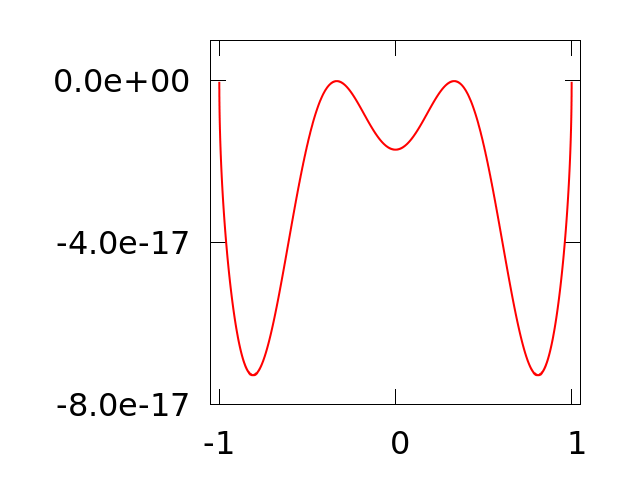}
}
\begin{picture}(-110,0)(0,0)
\put(-205,32){{\Large $=$}}
\put(-99,30){{\Large $\times$}}
\end{picture}
\caption{The factor $Q$ in \pRef{def_Q} in Salzer's case with $n + 1 = 100$ nodes
and $\wfc{f}{x} = \wfc{\cos}{10\, x}$. Notice that $Q$ is $\wfc{O}{\epsilon}$
and does not oscillate much.}
\label{figure_error_pol}
\end{figure}
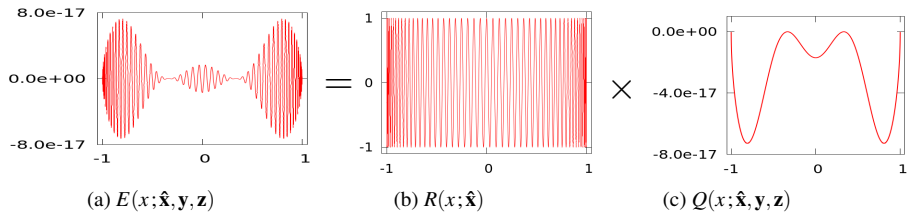

We  end this overview emphasizing that the Lipschitz constant plays
an important role in the accuracy of Step II for the second formula.
In Salzer's case, the function
\[
\wpr{x} \, :=  \, \frac{2^{n -1}}{\sqrt{1 - x^2}} \prod_{k = 0}^n \wlr{x - \hat{x}_k}
\]
is highly oscillating, with maximum absolute value close
to $1$, because $\wfc{R}{x;\wvec{x}^c} = - \wfc{\sin}{n \arccos x}$.
The function $\wfc{R}{x;\wvec{\hat{x}}}$ has simple zeros in $(-1,1)$
 and all of them are zeros of $\wpe{x}{y}$.
It follows that
\begin{equation}
\label{def_Q}
\wpq{x}{y}{z}  \, :=  \, \frac{\wpe{x}{y}}{\wpr{x}}
\end{equation}
is also a smooth function. This leads to the decomposition
\[
\wpe{x}{y} = \wpr{x} \times \wpq{x}{y}{z}.
\]
The factor $\wpq{x}{y}{z}$ depends on the interpolated function $f$.
Figure \ref{figure_error_pol}
illustrates graphically the decomposition $E = R \times Q$ for $\wfc{f}{x} = \wfc{\cos}{10 \, x}$
in Salzer's case

\begin{figure}[!h]
\subfloat[$\wpq{x}{y}{z}$ for $\wfc{\cos}{10^2 x}$]{
\includegraphics[viewport= 35 10 600 380, width=3.7cm, height=2.3cm]{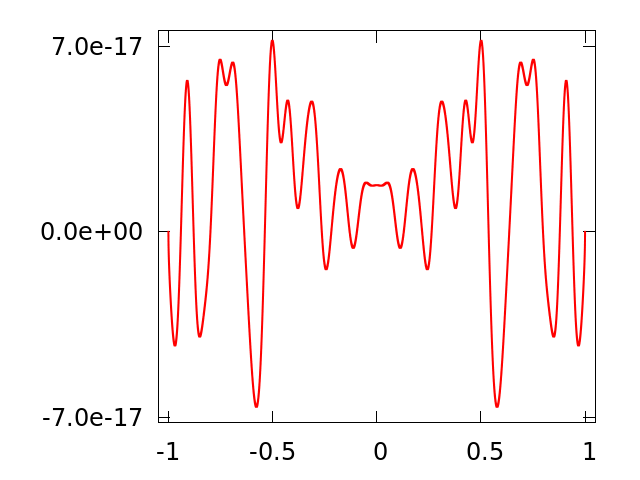}
}
\subfloat[$\wpq{x}{y}{z}$ for $\wfc{\cos}{10^3 x}$]{
\includegraphics[viewport= 20 10 600 380, width=3.7cm, height=2.3cm]{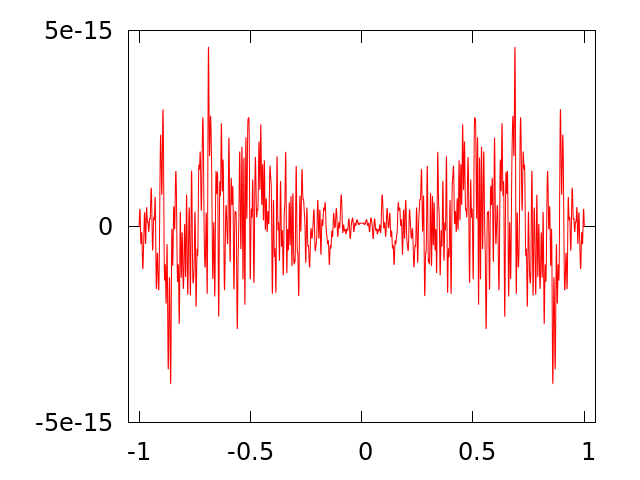}
}
\subfloat[$\wpq{x}{y}{z}$ for $\wfc{\cos}{10^4 x}$]{
\includegraphics[viewport= 20 10 560 380, width=3.7cm, height=2.3cm]{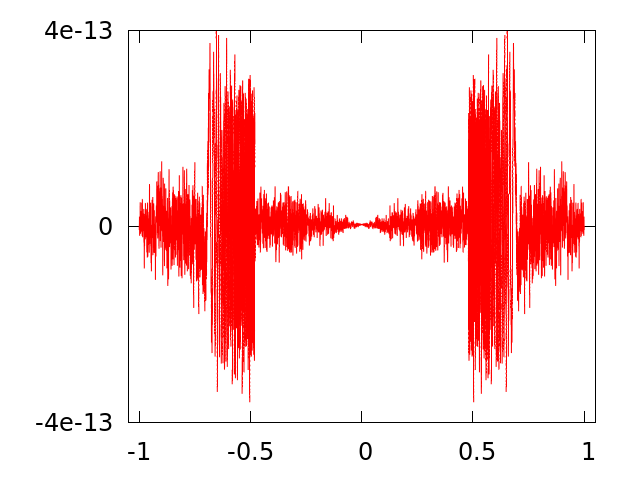}
}
\caption{The effect of high derivatives on the factor $Q$, for one million nodes.}
\label{figure_lipschitz_error}
\end{figure}

Figure \ref{figure_lipschitz_error} shows that larger Lipschitz constants
lead to larger amplitudes as well as larger frequencies for the errors. In the extreme
case given by the Lagrange polynomials with nodes $\wvec{\hat{x}}$, the
Lipschitz constant is of order $n^2$, and the bounds presented in this article
are not encouraging. In fact, \cite{MascCam} shows that the maximum backward error
for the second formula is of order $\epsilon n^2$ for Lagrange polynomials
in Salzer's case.

\section{Notation and conventions}
\label{section_notation}
Throughout the article, we consider intervals $[\hat{x}^-,\hat{x}^+]$, $[x^-,x^+]$
nodes $\wvec{\hat{x}}$ and $\wvec{x}$ and weights $\wvec{\hat{w}}$ and $\wvec{w}$, and readers
may find it convenient to have a copy of the next equations at hand while they read our arguments.
We convention the following:
\begin{eqnarray}
\label{first_cond}
x_k < x_{k+1}, \hspace{0.5cm} \hat{w}_k \neq 0 \hspace{0.5cm} \wrm{and} \hspace{0.5cm} w_k \neq 0, \hspace{0.5cm} \\
\nonumber
x_k \in \wlr{x^-,x^+} \  \wrm{if \ and \ only \ if } \ \hat{x}_k  \in \wlr{\hat{x}^-,\hat{x}^+}, \hspace{0.5cm} \\
x_k = x^- \ \  \wrm{if \ and \ only \ if } \ \hat{x}_k = \hat{x}^-, \hspace{1cm}
x_k = x^+ \ \ \wrm{if \ and \ only \ if } \ \hat{x}_k = \hat{x}^+, \hspace{0.5cm} \\
k^- \ \wrm{is \ the \ smallest \ } k \ \wrm{such \ that } \ x_k > x^-, \hspace{0.5cm} \\
k^+ \ \wrm{is \ the \ largest \ } k \ \wrm{such \ that } \ x_k < x^+ \hspace{0.5cm} \wrm{and}
\hspace{0.5cm} k^- \leq k^+.\hspace{0.5cm}
\label{last_cond}
\end{eqnarray}
We measure the difference  between $\wvec{x}$ and $\wvec{\hat{x}}$ in terms of
\begin{eqnarray}
\label{def_delta}
\delta_{kk} & := & \wfc{\delta_{kk}}{\wvec{\hat{x}},\wvec{x}} := 0 \hspace{1cm} \wrm{and} \hspace{1cm}
\delta_{jk} := \wfc{\delta_{jk}}{\wvec{\hat{x}},\wvec{x}} := \frac{x_j - x_k}{\hat{x}_j - \hat{x}_k} - 1, \\
\nonumber
\delta^-_j & := & \wfc{\delta^-_j}{\hat{x}^-, \wvec{\hat{x}}, x^-, \wvec{x}} :=
\frac{x^- - x_j}{\hat{x}^- - \hat{x}_j} - 1, \hspace{0.5cm}
\delta^+_j := \wfc{\delta^+_j}{\wvec{\hat{x}}, \hat{x}^+, \wvec{x}, x^+} :=
\frac{x^+ - x_j}{\hat{x}^+ - \hat{x}_j} - 1,
\end{eqnarray}
with the convention that $\delta^-_j :=  0$ and $\delta^+_j := 0$
in the exceptional cases $\hat{x}_j = \hat{x}^-$ and $\hat{x}_j = \hat{x}^+$.
We combine the $\delta^-_j$, $\delta_{jk}$ and $\delta^+_j$ in
\begin{equation}
\label{max_delta}
\delta := \max_{0 \leq j,k \leq n} \wset{\wabs{\delta^-_j}, \wabs{\delta_{jk}},
                                         \wabs{\delta^+_j}}.
\end{equation}
We also measure the errors in the nodes by
\begin{eqnarray}
\nonumber
\Delta^- & := & \sum_{j = 0}^n \max \wset{\wabs{\delta^-_j},\wabs{\delta_{j\wlr{k^-}}}},
\hspace{0.7cm}
\Delta_k :=  \sum_{j \neq k} \max \wset{\wabs{\delta_{jk}},\wabs{\delta_{j\wlr{k+1}}}}, \hspace{0.2cm} \\
\label{big_delta}
\Delta^+  & := & \sum_{j = 0}^n \max \wset{\wabs{\delta_{j\wlr{k^+}}}, \wabs{\delta_j^+}}
\hspace{0.7cm}
\Delta := \max_{k^- \, \leq k \, \leq \, \min\wset{k^+,n-1}} \wset{\Delta^-, \ \Delta_k, \ \Delta^+ }
\hspace{0.2cm}
\end{eqnarray}
($\Delta_k$ is defined for $0 \leq k < n$.)
The differences between the weights $\hat{\wvec{w}}$ used in computation and the weights
$\wfc{\lambda_k}{\hat{\wvec{x}}}$ corresponding to the rounded nodes according to
\pRef{def_lambda} are measured by
\begin{equation}
\label{def_zk}
z_k := \wfc{z_k}{\hat{\wvec{x}},\hat{\wvec{w}}} :=
\frac{\hat{w}_k - \wfc{\lambda_k}{\hat{\wvec{x}}}}{\wfc{\lambda_k}{\hat{\wvec{x}}}}
\hspace{1cm} \wrm{and} \hspace{1cm}
\zeta_k := \wfc{\zeta_k}{\wvec{w},\hat{\wvec{w}}} = \frac{w_k - \hat{w}_k}{\hat{w}_k},
\end{equation}
in which usually $w_k = \wfc{\lambda_k}{\hat{\wvec{x}}}$.
The Lebesgue constant is defined as
\begin{equation}
\label{def_lebesgue}
\Lambda := \Lambda_{x^-,x^+,\wvec{x}} := \sup_{x \in [x^-,x^+], \ \ \wvec{y} \neq 0} \frac{\wabs{\wfc{P}{x;\wvec{x},\wvec{y}}}}{\wnorm{\wvec{y}}_\infty},
\end{equation}
where $\wfc{P}{x;\wvec{x},\wvec{y}} \ \wrm{is \ the \ } n \wrm{th \ degree \ polynomial \ that \ interpolates \ } \wvec{y} \ \wrm{at} \ \wvec{x}$,
Note that this definition implies that if $Q$ is a $n$th degree polynomial, then
\begin{equation}
\label{lebesgue_prop}
\wabs{\wfc{Q}{x_k}} \leq M \ \ \wrm{for} \ \ 0 \leq k \leq n \ \ \Rightarrow \ \ \wabs{\wfc{Q}{x}}
\leq \Lambda_{x^-,x^+,\wvec{x}} M \ \ \wrm{for} \ \ x \in [x^+,x^-],
\end{equation}
because $\wfc{Q}{x} =  \wfc{P}{x,\wvec{x},\wfc{Q}{\wvec{x}}}$.
The $k$th Lagrange polynomial with nodes $\wvec{x}$ is given by
\begin{equation}
\label{def_lagrange}
\wplagrx{k}{x} := \wfc{\lambda_k}{\wvec{x}} \prod_{j \neq k} \wlr{x - x_k}.
\end{equation}

The $z_k$, $\zeta_k$ and $\delta_{jk}$ are related by the following lemma:
\pbLemma{lemBoundZk}
If $a_k := \sum_{j \neq k} \wabs{\delta_{jk}} < 1$ then $z_k = \wfc{z_k}{\hat{\wvec{x}}, \wfc{\lambda_k}{\wvec{x}}}$ and
$\zeta_k = \wfc{\zeta_k}{\wfc{\lambda_k}{\hat{\wvec{x}}}, \wfc{\lambda_k}{\wvec{x}}}$ satisfy
\pbTClaim{boundZk}
0  \leq \ z_k + \sum_{j \neq k} \delta_{jk} \leq \frac{a_k^2}{1 - a_k}
\pStop{,} \hspace{0.8cm}
\wabs{z_k} \leq \frac{a_k}{1 - a_k}
\hspace{0.8cm} \pStop{and} \hspace{0.8cm}
\wabs{\zeta_k} \leq \frac{a_k}{1 - a_k}.
\peTClaim{boundZk}
\peLemma{lemBoundZk}

Using $\Delta$, we can bound products of the form
\[
\prod \wlr{x - x_k}/\wlr{\hat{x} - \hat{x}_k}:
\]
\pbLemma{lemBackwardProduct}
Under the conditions \pRef{first_cond}--\pRef{last_cond}, if $\Delta < 1$ then
for every $\hat{x} \in [\hat{x}^-,\hat{x}^+]$  there exists $\wfc{x}{\hat{x}}  \in[x^-,x^+]$ such that
\begin{equation}
\label{delta_xu}
\wabs{\wfc{x}{\hat{x}}  - \hat{x}} \leq
\max \wset{ \wnorm{\wvec{x} - \wvec{\hat{x}}}_\infty, \ \wabs{x^- - \hat{x}^-},  \ \wabs{x^+ - \hat{x}^+}}
\end{equation}
and, for the same $\wfc{x}{\hat{x}}$ and  every $K \subset \wset{0,1,\dots,n}$, there exists
$\beta_K$ such that
\begin{equation}
\label{backward_product}
\wabs{\beta_K} \leq \frac{\Delta}{1 - \Delta} \hspace{1cm} \wrm{and} \hspace{1cm}
\prod_{k \in K} \wlr{\hat{x} - \hat{x}_k} = \wlr{1 + \beta_K} \prod_{k \in K} \wlr{\wfc{x}{\hat{x}}  - x_k}.
\end{equation}
\peLemma{lemBackwardProduct}

In Salzer's case, we can bound $\delta$, $\Delta$, $\wvec{z}$  and $\bm{\zeta}$ in terms of $n$ and the rounding errors in the nodes,
and show that our theory applies even to $n$ in the million range:


\pbLemma{lemSalzersConstants}
In \hyperlink{link_salzer}{Salzer's case} we have the bounds
in Table \ref{table_salzers_bounds}, where $\wvec{z} = \wfc{\wvec{z}}{\hat{\wvec{x}},\wfc{\lambda}{\wvec{x}}}$
and $\bm{\zeta} = \wfc{\bm{\zeta}}{\wfc{\lambda}{\hat{\wvec{x}}},\wfc{\lambda}{\wvec{x}}}$:
\begin{table}[!h]
\caption{Upper bounds for Salzer's case}
\centering
\begin{tabular}{c|llc}
\hline\\[-0.27cm]
         &   &       & absolute\\
         &   &  upper bound for $n$     &upper bound\\
\hline\\[-0.25cm]
$\delta$ & $0.40897 \wnorm{\wvec{\hat{x}}^c - \wvec{x}^c}_\infty n^2$          & $1.8813 \times 10^{-16}\, n^2$  &$7.5252 \times 10^{-4}$ \\[0.05cm]
$\Delta$ & $2.7267  \wnorm{\wvec{\hat{x}}^c - \wvec{x}^c}_\infty n^2$          & $1.2543 \times 10^{-15}\, n^2$  &$5.0172 \times 10^{-3}$ \\[0.05cm]
$\max \sum_{j} \wabs{\delta_{jk}}$& $2.4502 \wnorm{\wvec{\hat{x}}^c - \wvec{x}^c}_\infty n^2$ & $1.1271 \times 10^{-15}\, n^2$  &$4.5084 \times 10^{-3}$ \\[0.05cm]
$\wnorm{\wvec{z}}_\infty$       & $2.4624 \wnorm{\wvec{\hat{x}}^c - \wvec{x}^c}_\infty n^2$          & $1.1328 \times 10^{-15}\, n^2$  &$4.5312 \times 10^{-3}$ \\[0.05cm]
$\wnorm{\bm{\zeta}}_\infty$       & $2.4624 \wnorm{\wvec{\hat{x}}^c - \wvec{x}^c}_\infty n^2$          & $1.1328 \times 10^{-15}\, n^2$  &$4.5312 \times 10^{-3}$ \\[0.05cm]
$\wnorm{\wvec{z}}_1$            & $1.4219 \wnorm{\wvec{\hat{x}}^c - \wvec{x}^c}_\infty n^2 \wlr{3 + \log n}^2$ &
                                  $6.5406 \times 10^{-16} \, n^2 \wlr{3 + \log n}^2$&
                                  $8.0202 \times 10^{-1} $ \\[0.05cm]
\end{tabular}
\label{table_salzers_bounds}
\end{table}
\peLemma{lemSalzersConstants}

\section{Perturbations in the nodes of polynomials}
\label{section_polynomials}
There are at least four reasonable concepts of ``the interpolating polynomial of $f$'' when
we take into account rounding errors in the nodes:
\[
\wfc{P}{x;\wvec{x},\wfc{f}{\wvec{x}}}, \hspace{1cm}
\wfc{P}{x;\wvec{x},\wfc{f}{\hat{\wvec{x}}}}, \hspace{1cm}
\wfc{P}{x;\wvec{\hat{x}}, \wfc{f}{\wvec{x}}}
\hspace{1cm} \wrm{and} \hspace{1cm}
\wfc{P}{x;\wvec{\hat{x}},\wfc{f}{\wvec{\hat{x}}}}.
\]
The difference $\wfc{f}{\wvec{\hat{x}}} - \wfc{f}{\wvec{x}}$ is of order
$L \wsupdx$, where $L$ is $f$'s Lipschitz constant.
This leads to a difference
$\wfc{P}{x;\wvec{x},\wfc{f}{\wvec{x}}} - \wfc{P}{x;\wvec{x},\wfc{f}{\wvec{\hat{x}}}}$
of order $\Lambda_{x^-,x^+,\wvec{x}} L \wsupdx$, and this
section shows that when $f$ is well approximated by polynomials this is the overall
order of magnitude of the difference $\wfc{P}{x;\wvec{x},\wvec{y}} - \wfc{P}{x;\hat{\wvec{x}},\wvec{y}}$,
as suggested in the vertical edge in Figure \ref{figure_bounds}.
We also show that the usual error estimate for Lagrange interpolation, namely,
\begin{equation}
\label{lagrange_bound}
\frac{\wfc{f^{(n+1)}}{\xi}}{\wlr{n + 1}!} \prod_{k = 0}^n \wlr{x - x_k},
\end{equation}
does not change much when we replace $\wvec{x}$ by $\hat{\wvec{x}}$ and the $\Delta$ in \pRef{big_delta} is small.
Therefore, $\wfc{P}{x;\wvec{\hat{x}},\wfc{f}{\wvec{\hat{x}}}}$ is as good an approximation
of $f$ as $\wfc{P}{x;\wvec{x},\wfc{f}{\wvec{x}}}$ if we
consider \pRef{lagrange_bound} as a measure of the degree of approximation.
Moreover, we show that when $\Delta$ is small the Lebesgue constant
with respect to the nodes $\wvec{\hat{x}}$ is roughly the same as the Lebesgue
constant corresponding to the nodes $\wvec{x}$.
The overall conclusion of this section is that it is reasonable to consider the differences
\[
\wfc{p}{x;\wvec{\hat{x}},\wvec{y},\hat{\wvec{w}}} - \wfc{p}{x;\wvec{\hat{x}},\wvec{y},\wfc{\lambda}{\wvec{\hat{x}}}}
\hspace{0.5cm} \wrm{and} \hspace{0.5cm}
\wfc{q}{x;\wvec{\hat{x}},\wvec{y},\hat{\wvec{w}}} - \wfc{q}{x;\wvec{\hat{x}},\wvec{y},\wfc{\lambda}{\wvec{\hat{x}}}}
\]
as measures of the errors in Step II for the first and second formula,
as in Salzer's case.
The informal statements above are formalized by the following results:
\pbLemma{lemLebesguePerturbation}
Under the conditions \pRef{first_cond}--\pRef{last_cond},
if  $\wvec{z} = \wfc{\wvec{z}}{\hat{\wvec{x}}, \wfc{\lambda}{\wvec{x}}}$ satisfies
\begin{equation}
\label{hypo_lebesgue_perturbation}
\delta <  \frac{1 - \wnorm{\wvec{z}}_\infty}{\Lambda_{x^-,x^+,\wvec{x}} } - \wnorm{\wvec{z}}_\infty
\end{equation}
then
\begin{equation}
\label{bound_lebesgue}
\Lambda_{\hat{x}^-,\hat{x}^+,\wvec{\hat{x}}} \leq \frac{1 + \delta}{1 - \wnorm{\wvec{z}}_\infty -
\wlr{\delta + \wnorm{\wvec{z}}_\infty} \Lambda_{x^-,x^+,\wvec{x}}} \Lambda_{x^-,x^+,\wvec{x}}.
\end{equation}
In particular, in Salzer's case with $n \leq 10^6$ we have
\begin{equation}
\label{bound_lebesgue_salzer}
\Lambda_{-1,1,\wvec{\hat{x}}^c} \leq 1.0629 \Lambda_{-1,1,\wvec{x}^c} \leq 0.67667 \log n + 1.0236
\hspace{0.35cm} \wrm{and} \hspace{0.35cm}
\Lambda_{-1,1,\wvec{\hat{x}}^c} \leq 10.841.
\end{equation}
\peLemma{lemLebesguePerturbation}
\pbLemma{lemWarmUp}
Under the conditions \pRef{first_cond}--\pRef{last_cond}, if $\Delta < 1$ then
the usual estimate \pRef{lagrange_bound} for $[\hat{x}^-,\hat{x}^+]$ and
$\wvec{\hat{x}}$ is not much larger than the same estimate for $[x^-,x^+]$
and $\wvec{x}$, namely,
\begin{equation}
\label{bound_warm_up}
\max_{\hat{x} \in [\hat{x}^-,\hat{x}^+]} \wabs{\prod_{k = 0}^n \wlr{\hat{x} - \hat{x}_k}} \leq
\frac{1}{1 - \Delta} \max_{x \in [x^-,x^+]} \wabs{\prod_{k = 0}^n \wlr{x - x_k}}.
\end{equation}
In Salzer's case,
\begin{equation}
\label{cheby_prod}
\max_{\hat{x} \in [-1,1]}
\wabs{\prod_{k = 0}^n \wlr{\hat{x} - \hat{x}^{(c)}_k}} \leq \wlr{1 + 1.2607 \times 10^{-15} n^2} \max_{x \in [-1,1]} \wabs{\prod_{k = 0}^n \wlr{x - x^{(c)}_k}}
\end{equation}
\[
\leq
2^{1 - n} \wlr{1 + 1.2607 \times 10^{-15} n^2} \leq 2.0101 \times 2^{- n}.
\]
\peLemma{lemWarmUp}

It is well known that the best possible value for the left hand
side of \pRef{cheby_prod}, among all sets of nodes, is at most $2^{-n}$ (see \cite{SALZER}.) Therefore,
the rounded nodes $\wvec{\hat{x}}^c$ are nearly optimal concerning the size of the product
in the bound \pRef{lagrange_bound}. More generally, Lemma \ref{lemWarmUp} shows
that, for reasonably rounded nodes, the accuracy for interpolation
of variations ($\wrm{a}$) and ($\wrm{c}$) in \pRef{eq_variations}
cannot be distinguished solely on basis of the traditional estimate \pRef{lagrange_bound}.
Finally, we present a theorem formalising the bound presented on the vertical edge in Figure \ref{figure_bounds}:

\pbTheorem{thmBestApprox}
Suppose that $\hat{x}^- \geq x^-$ and $\hat{x}^+ \leq x^+$
and consider a function $f: \wrone{} \rightarrow \wrone{}$ with Lipschitz constant $L$
and such that $y_k = \wfc{f}{x_k}$.
If $Q$ is a polynomial of degree at most $n$ and $x \in [\hat{x}^-,\hat{x}^+]$
then
\begin{equation}
\label{best_approx}
\wabs{\wfc{P}{x;\wvec{\hat{x}},\wvec{y}} - \wfc{P}{x;\wvec{x},\wvec{y}}} \leq
L \Lambda_{\hat{x}^-,\hat{x}^+,\wvec{\hat{x}}} \wsupdx +
M \wlr{\Lambda_{\hat{x}^-,\hat{x}^+,\wvec{x}} + \Lambda_{\hat{x}^-,\hat{x}^+,\wvec{\hat{x}}}} ,
\end{equation}
where
\begin{equation}
\label{def_q}
M := \max_{x \in \wset{x_0,x_1,\dots, x_n} \cup \wset{\hat{x}_0, \hat{x}_1, \dots, \hat{x}_n}} \wabs{\wfc{Q}{x} - \wfc{f}{x}}.
\end{equation}
In Salzer's case,  for all $x \in [-1,1]$, we have
\begin{equation}
\label{bound_salzer_approx}
\wabs{\wfc{P}{x;\wvec{\hat{x}},\wvec{y}} - \wfc{P}{x;\wvec{x},\wvec{y}}} \leq
\wlr{0.7 \log n + 1} L  \wsupdx +
M \wlr{2 + 1.4 \log n}.
\end{equation}
\peTheorem{thmBestApprox}

\section{Bounds for the first formula}
\label{section_first}
Here we discuss the forward and backward stability of the first formula in
Steps II and III. The first issue we address is the appropriate concept of backward
stability when we allow for rounded nodes.
We then present upper bounds
on the errors in Steps II and III, and lower bounds for Step II, for the first formula. In particular, we
show that, when errors in the nodes are taken into account,
the overall error in steps II and III in Salzer's case
is of order $\epsilon n^2$ for this formula.

We now explain that when the nodes are perturbed we cannot restrict
ourselves to perturbation in the function values in order to
prove backward stability.
In fact, let us consider the Lagrange polynomials,
which can be written in first barycentric form as
\[
\wfc{\ell_k}{x;\wvec{x}} = \wfc{p}{x, \wvec{x}, \wvec{e}^k y_k, \wfcc{\lambda}{\wvec{x}}},
\]
where $y_k = 1$ and $\wvec{e}^k \in \wrn{n+1}$ is the vector with
$\wvec{e}^{(k)}_k = 1$ and $\wvec{e}^{(k)}_j = 0$ for $j \neq k$.
If we were to consider only
relative perturbations in $\wvec{y}$, then given rounded
nodes $\wvec{\hat{x}}$ we would need to find $\hat{y}_k = \wlr{1 + \beta_k} y_k = 1 + \beta_k$
such that
\[
\wfc{p}{x, \wvec{\hat{x}}, \wvec{e}^k, \wfcc{\lambda}{\wvec{x}}} =
\wfc{p}{x, \wvec{x}, \wvec{e}^k \wlr{1 + \beta_k}, \wfcc{\lambda}{\wvec{x}}}.
\]
This equation leads to
\[
\beta_k = \prod_{i \neq k} \frac{x - \hat{x}_k}{x - x_k} - 1,
\]
and given an arbitrarily small $\delta > 0$, we could take
$k = 1$, $\hat{x}_j = x_j$ for
$k \neq 1$, $\hat{x}_1 = x_1 + \delta$ and obtain
\[
\beta_1 = \frac{x - x_1 - \delta}{x - x_1} - 1 = \frac{-\delta}{x - x_1}.
\]
We could then make $\beta_1$ arbitrarily large by taking $x$ close enough to
$x_1$. Therefore, we cannot build a backward stability theory for
the first formula with perturbed nodes relying only on perturbations
of the function values.
Fortunately, the next theorem shows that we can get meaningful results
if we allow for perturbations in $x$:

\pbTheorem{thmFirstBack}
Under the conditions \pRef{first_cond}--\pRef{last_cond}, if $\Delta < 1$
and the machine precision $\epsilon$ is such that $\wlr{3 n + 5} \epsilon < 1$,
then for every $\hat{x} \in [\hat{x}^-,\hat{x}^+]$ and
$\wvec{y} \in \wrn{n+1}$ there exists
$x \in [x^-,x^+]$ and $\bm{\beta}, \bm{\nu} \in \wrn{n+1}$ such that
\[
\wnorm{\bm{\beta}}_\infty \leq \Delta / \wlr{1 - \Delta}, \hspace{2cm}
\wnorm{\bm{\nu}}_\infty \leq \frac{\wlr{3 n + 5} \epsilon}{1 - \wlr{3 n + 5} \epsilon}
\]
and
\begin{equation}
\label{dx_first_back}
\wabs{\hat{x} - x} \leq \max \wset{\wsupdx,
\ \wabs{\hat{x}^- - x^-}, \ \ \wabs{\hat{x}^+ - x^+} },
\end{equation}
and the vector $\tilde{\wvec{y}}$ with
\begin{equation}
\label{y_tilde}
\tilde{y}_k := y_k  \wlr{1 + \beta_k}\wlr{1 + \nu_k}\wlr{1 + \frac{\hat{w}_k - w_k}{w_k}}
\end{equation}
satisfies
\begin{equation}
\label{first_back_match}
\wrounde{\wfc{p}{\hat{x}, \wvec{\hat{x}}, \wvec{y}, \hat{\wvec{w}}}} =
\wfc{p}{x, \wvec{x},\tilde{\wvec{y}}, \wvec{w}}.
\end{equation}
\peTheorem{thmFirstBack}
In words, Theorem \ref{thmFirstBack} shows that the first formula is
backward stable in Steps II and III, in the broader sense which allows 
also for perturbations in $x$.

Our analysis of the forward stability of the first formula in Steps II and III is more complete
than the analysis of the backward stability, because it also yields an estimate on the error,
which is given in equation \pRef{lower_bound_first} in the following theorem:
\pbTheorem{thmFirstForward}
If $\hat{x}\in [\hat{x}^-,\hat{x}^+]$ and the machine precision $\epsilon$ is such that $\wlr{3 n + 5} \epsilon < 1$,
then the first formula $p$ in \pRef{first_formula} satisfies
\begin{equation}
\label{bound_first}
\wabs{\wrounde{\wfc{p}{\hat{x},\wvec{\hat{x}},\wvec{y},\hat{\wvec{w}}}} - \wpit{\hat{x}}{y}} \leq
\Lambda_{\hat{x}^-,\hat{x}^+,\wvec{\hat{x}}} \max_{0 \leq j \leq n} \wabs{y_j \, \frac{z_j + \wlr{3 n + 5} \epsilon}{1 - \wlr{3 n + 5} \epsilon}},
\end{equation}
for $\wvec{z} = \wfc{\wvec{z}}{\hat{\wvec{x}}, \hat{\wvec{w}}}$.
Moreover, for $\hat{x}, \hat{x}_k \in [\hat{x}^-,\hat{x}^+]$,
there exists a vector $\wvec{\nu} \in \wrn{n+1}$ with
\begin{equation}
\label{bound_nuk}
\wabs{z_j} \wlr{1 - \wlr{3 n + 5} \epsilon} - \frac{\wlr{3 n + 5} \epsilon}{1 - \wlr{3 n + 5} \epsilon}
\leq \wabs{\nu_j} \leq
\frac{\wabs{z_j} + \wlr{3 n + 5} \epsilon}{1 - \wlr{3 n + 5} \epsilon}
\end{equation}
such that, for $\wvec{\hat{x}}' = \wset{\hat{x}_0,\dots,\hat{x}_n} \setminus \wset{\hat{x}_k}$, we have
\begin{equation}
\label{lower_bound_first}
\wabs{\wlr{\wrounde{\wfc{p}{\hat{x},\wvec{\hat{x}},\wvec{y},\hat{\wvec{w}}}} - \wpit{\hat{x}}{y}}  - y_k \nu_k}
\leq
\Lambda_{\hat{x}^-,\hat{x}^+,\wvec{\hat{x}'}} \wabs{\hat{x} - \hat{x}_k}
\max_{j \neq k} \wabs{ \frac{y_j \nu_j - y_k \nu_k}{\hat{x}_j - \hat{x}_k}}.
\end{equation}
\peTheorem{thmFirstForward}

In Salzer's case, Lemma \ref{lemSalzersConstants} shows that
the upper bound \pRef{bound_first} is of order $\epsilon n^2 \log n \wabs{y_j}$.
Of course, in itself, this upper bound does not imply that the backward errors will be
of order $\epsilon n^2 \log n \wnorm{\wvec{y}}_\infty$ in this case. In fact, in the usual situations
in which $x$ is not very close to the nodes, the backward error will be much smaller than the
right hand side of \pRef{bound_first}.
However, Table 1 in \cite{MascCam}
provides strong empirical evidence that $\wnorm{\wvec{z}}_\infty$ is at least $0.01 \epsilon n^2$ when we round the Chebyshev
points as usual (note that $\wabs{\zeta_k} = \wabs{z_k/\wlr{1 + z_k}} \approx \wabs{z_k}$.)
In this case, equation \pRef{bound_nuk} bounds $\nu_j$ below by an estimate of order $\epsilon n^2$ and
equation \pRef{lower_bound_first} suggests that whenever the $y_k$ corresponding to
$\wabs{z_k} \approx \wnorm{\wvec{z}}_\infty$ is not small,
it is likely that we will incur in errors of magnitude $\epsilon n^2 \wnorm{\wvec{y}}_\infty$ when we evaluate
the first barycentric formula for $\hat{x}$ very close to $\hat{x}_k$.
For instance, when $n = 10^6$ and $\epsilon = 2.3 \times 10^{-16}$ we have
$0.01 \epsilon n^2 =  10^{-2} \times 2.3 \times 10^{-16} \times 10^{12} = 2.3 \times 10^{-6}$,
and this value agrees remarkably well with the maximum error of $6 \times 10^{-6}$ for the
sine function presented in the Table 5 of  \cite{Masc}, with the data displayed in Figure \ref{figure_least_squares}, and also with
the maximum error of $2 \times 10^{-6}$ for the function $\wfc{\cos}{100x}$
presented in Tables \ref{table_stepII_error} and \ref{table_overall_error}
in Subsection \ref{subsection_experiments}.

\section{Bounds for the second formula}
\label{section_second}
In this section we bound the errors in Step II and Step III  in Figure \ref{figure_steps} for the
second barycentric formula \pRef{second_formula} in terms of the Error Polynomial
in \pRef{def_error_pol}. The Error Polynomial is similar to a function presented by
Werner \cite{WERNER}, which expresses its results in terms of divided differences.
Our approach is slightly different: we use the Error Polynomial in combination with our bounds on
$\wvec{z}$
in order to have a unified picture for both barycentric formulae and to obtain more explicit bounds.
The next theorem relates the Error Polynomial,
the Lebesgue constant $\Lambda_{\hat{x}^-,\hat{x}^+,\wvec{\hat{x}}}$,
and the forward error in Step II for the second formula:
\pbTheorem{thmMain}
Consider $\wvec{z} = \wfc{\wvec{z}}{\hat{\wvec{x}},\hat{\wvec{w}}}$.
If $\hat{x} \in [\hat{x}^-,\hat{x}^+]$ and
$\Lambda_{\hat{x}^-,\hat{x}^+,\wvec{\hat{x}}} \wnorm{\wfc{\wvec{z}}{\hat{\wvec{x}},\hat{\wvec{w}}}}_\infty < 1$ then
the second formula $q$ in \pRef{second_formula} satisfies
\begin{equation}
\label{good_news}
\frac{\wabs{ \wpe{\hat{x}}{y}}}
{{1 + \Lambda_{\hat{x}^-,\hat{x}^+,\wvec{\hat{x}}} \wnorm{\wvec{z}}_\infty}} \ \leq \
\ \wabs{\wfc{q}{\hat{x},\wvec{\hat{x}},\wvec{y},\hat{\wvec{w}}} - \wpit{\hat{x}}{y} }
\ \leq \ \frac{\wabs{ \wpe{\hat{x}}{y}}}
{{1 - \Lambda_{\hat{x}^-,\hat{x}^+,\wvec{\hat{x}}} \wnorm{\wvec{z}}_\infty}}.
\end{equation}
\peTheorem{thmMain}
\noindent
We now present an empirical way to bound $\wabs{ \wpe{\hat{x}}{y}}$, and after
this empirical bound we present a theoretical bound of order $L \epsilon n \log^2 n$ on the forward error of
the second barycentric formula.
Given $\wvec{y}$ and $\wvec{z}$, we can compute bounds on the Error Polynomial in terms of its values at
convenient points $c_k$. In fact, given points $\wset{c_0,c_1, \dots c_m}$
we define
\begin{equation}
\label{vk}
C := \wset{\hat{x}_k, \ 0 \leq k \leq n} \  \bigcup \ \wset{c_k, \ 1 \leq k \leq m}
\end{equation}
and the identity $\wpe{\hat{x}_k}{y} = 0$ and the definition of Lebesgue constant \pRef{def_lebesgue} yield
\begin{equation}
\label{bound_2n}
\max_{\hat{x} \in [\hat{x}^-, \hat{x}^+]} \wabs{\wpe{\hat{x}}{y}} \leq
\Lambda_{\hat{x}^-,\hat{x}^+,C} \max_{1 \leq k \leq m} \wabs{\wpe{c_k}{y}}
\end{equation}
when $m \geq n$. The right-hand side of \pRef{bound_2n} overestimates the left-hand side by at most
a factor of $\Lambda_{\hat{x}^-,\hat{x}^+,C}$. In Salzer's case, if we choose
$C$ as the rounded Chebyshev nodes corresponding to $2n$ then Lemma
\ref{lemLebesguePerturbation} shows that $\Lambda_{\hat{x}^-,\hat{x}^+,C} \leq 0.7 \log (2n) + 1.1$
for $10 \leq n \leq 1.000.000$ and the right hand side \pRef{bound_2n}
overestimates the left hand side by a $\wfc{O}{\log n}$ factor.

Equation \pRef{bound_2n} leads to computable bounds for the forward errors for
relevant classes of functions.
We can bound the error in Step II for a Lipschitz function $f$
by evaluating $\wvec{z}$
and maximizing the linear functions of $\wvec{y}$ given by
\begin{equation}
\label{def_hk}
\wfc{h_k}{\wvec{y}} := \wpe{c_k}{y} = \sum_{j = 0}^n a_{k,j} y_j,
\end{equation}
where the $a_{k,j}$ are defined using $\wplagr{j}{x}$, the $j$-th Lagrange polynomial with nodes
$\wvec{\hat{x}}$,
\begin{equation}
\label{def_akj}
a_{k,j} := \wfc{a_{k,j}}{C, \wvec{z}} :=
\wlr{z_j  - \wpit{c_k}{\wvec{z}}} \wplagr{j}{c_k},
\end{equation}
subject to linear constraints of the form
\begin{equation}
\label{constraint_y}
y_j - y_{j-1} \leq M \wlr{ \hat{x}_j - \hat{x}_{j - 1} } \hspace{1cm} \wrm{and} \hspace{1cm}
y_{j-1} - y_j \leq M \wlr{ \hat{x}_j - \hat{x}_{j - 1} },
\end{equation}
where we can choose the constant $M$ as $f$'s Lipschitz constant.
The linear programming problem with objective function \pRef{def_hk} and constraints
\pRef{constraint_y} can be solved in closed form. Its solution leads to the following lemma:

\pbLemma{lemLp}
Consider $a_{k,j}$ in \pRef{def_akj} and  $M \in \wrone{}$.
If $\wvec{y}$ satisfies \pRef{constraint_y} for
$1 \leq j \leq n$ then
\[
\wabs{\wpe{c_k}{y}} \leq
M \sum_{i=1}^n \wlr{\hat{x}_{i} - \hat{x}_{i-1}} \wabs{\sum_{j = i}^n a_{k,j}} .
\]
\peLemma{lemLp}
This lemma with
\[
M =  \max_{1 \leq k \leq n} \wabs{\frac{y_k - y_{k-1}}{\hat{x}_k - \hat{x}_{k-1}}}
\]
 and equation \pRef{bound_2n} lead to
\begin{equation}
\label{lip_bound}
\max_{x \in [\hat{x}^-,\hat{x}^+]} \wabs{\wpe{x}{y}} \leq
\Lambda_{\hat{x}^-,\hat{x}^+,C} \wlr{ \max_{1 \leq k \leq n} \wabs{\frac{y_k - y_{k-1}}{\hat{x}_k - \hat{x}_{k-1}}} }
\wfc{b_n}{C,\wvec{z}},
\end{equation}
for $C$ in \pRef{vk} and
\begin{equation}
\label{def_bn}
\wfc{b_n}{C,\wvec{z}} := \max_{1 \leq k \leq n} \sum_{i=1}^n \wlr{\hat{x}_{i} - \hat{x}_{i-1}} \wabs{\sum_{j = i}^n \wfc{a_{k,j}}{C,\wvec{z}}}.
\end{equation}
When $C$ has $\wfc{O}{n}$ elements we can evaluate $\wfc{b_n}{C,\wvec{z}}$ in $\wfc{O}{n^2}$ operations.
 In fact, we can compute all the weights $\wfcc{\lambda_k}{\wvec{\hat{x}}}$ and $\wvec{z}$ in $\wfc{O}{n^2}$ operations.
 Next, we can evaluate $\wpit{c_k}{z}$ and $\wplagr{k}{c_k}$,
 before obtaining all the $a_{k,j}$ in \pRef{def_akj} for a given $k$ in $\wfc{O}{n}$ operations using the identity
\begin{equation}
\label{lagrange_jk}
\wplagr{j}{x} =
\frac{\wfcc{\lambda_j}{\wvec{\hat{x}}}}{\wfcc{\lambda_k}{\wvec{\hat{x}}}} \frac{x - \hat{x}_k}{x - \hat{x}_j}
\wplagr{k}{x}.
\end{equation}
For each $n$, set $C$ and vector $\wvec{z}$, we obtain
 a single number $\wfc{b_n}{C,\wvec{z}}$, which we can compute off-line.
 In Salzer's case, our computations with $C$ formed by the rounded Chebyshev
 points corresponding to $2n$ lead to the $\wfc{b_n}{C,\wvec{z}}$ in Table \ref{table_bn}.

 %
 %

 \begin{table}[!h]
\caption{$\wfc{b_n}{C,\wvec{z}}$ in Salzer's case with
$C$ formed by the rounded Chebyshev points corresponding to $2n$}
\centering
\begin{tabular}{c|c|c|c| c|c|c|c|c|c}
\hline\\[-0.27cm]
$n$ & $\wfc{b_n}{C,\wvec{z}}$ & $n$ & $\wfc{b_n}{C,\wvec{z}}$ &
$n$ & $\wfc{b_n}{C,\wvec{z}}$ & $n$ & $\wfc{b_n}{C,\wvec{z}}$ &
$n$ & $\wfc{b_n}{C,\wvec{z}}$  \\[-0.01cm]
\hline\\[-0.2cm]
  10 &  5.2e-17 &  100 &  2.1e-16 & 1.000 &  4.5e-16 & 10.000 &  5.5e-16 & 100.000 &  9.0e-16 \\
  20 &  1.2e-16 &  200 &  2.9e-16 & 2.000 &  5.0e-16 & 20.000 &  8.6e-16 & 200.000 &  8.7e-16 \\
  30 &  1.4e-16 &  300 &  3.6e-16 & 3.000 &  4.7e-16 & 30.000 &  7.1e-16 & 300.000 &  1.3e-15 \\
  40 &  1.1e-16 &  400 &  3.0e-16 & 4.000 &  5.3e-16 & 40.000 &  8.1e-16 & 400.000 &  1.1e-15 \\
  50 &  1.7e-16 &  500 &  3.4e-16 & 5.000 &  5.7e-16 & 50.000 &  8.4e-16 & 500.000 &  1.2e-15 \\
  60 &  1.7e-16 &  600 &  3.9e-16 & 6.000 &  6.2e-16 & 60.000 &  9.1e-16 & 600.000 &  1.2e-15 \\
  70 &  1.3e-16 &  700 &  3.9e-16 & 7.000 &  6.6e-16 & 70.000 &  7.7e-16 & 700.000 &  1.4e-15 \\
  80 &  1.6e-16 &  800 &  5.6e-16 & 8.000 &  7.1e-16 & 80.000 &  8.0e-16 & 800.000 &  1.5e-15 \\
  90 &  2.0e-16 &  900 &  6.3e-16 & 9.000 &  7.9e-16 & 90.000 &  8.9e-16 & 900.000 &  1.2e-15
\end{tabular}
\label{table_bn}
\end{table}

We now present a bound for the error in Step II for the second formula based only on the $z_k$.
This bound does not take the cancelation of rounding errors into account
and, consequently, it is worse than the bound obtained by combining Table \ref{table_bn}, equation \pRef{lip_bound} and Theorem \ref{thmMain}.
\pbTheorem{thmLipschitz}
Let $f: [\hat{x}^-,\hat{x}^+] \rightarrow \wrone{}$ be a function with Lipschitz constant $L$
such that $y_k = \wfc{f}{\hat{x}_k}$. If $\hat{x} \in [\hat{x}^-,\hat{x}^+]$ then
\begin{eqnarray}
\label{lipschitz_bound_a}
\wabs{\wpeu{\hat{x}}{y}} & \leq & L \, \wfc{\tau}{\wvec{\hat{x}}} \wnorm{\wvec{z}}_1 +
\Lambda_{\hat{x}^-,\hat{x}^+,\wvec{\hat{x}}} \ \wnorm{\wvec{z}}_\infty \, \sup_{s \in [\hat{x}^-,\hat{x}^+]} \wabs{\wfc{f}{s} - \wfc{P}{s,\wvec{\hat{x}},\wvec{y}}}\\
\label{lipschitz_bound_b}
&\leq& L \, \wlr{\wfc{\tau}{\wvec{\hat{x}}} \wnorm{\wvec{z}}_1 +
\Lambda_{\hat{x}^-,\hat{x}^+,\wvec{\hat{x}}} \wlr{1 + \Lambda_{\hat{x}^-,\hat{x}^+,\wvec{\hat{x}}}}
\wnorm{\wvec{z}}_\infty \frac{\wlr{\hat{x}^+ - \hat{x}^-} \pi}{4 \wlr{n + 1}}},
\end{eqnarray}
where
\begin{equation}
\label{def_tau}
\wfc{\tau}{\wvec{\hat{x}}} := \max_{\hat{x} \in [\hat{x}^-,\hat{x}^+]}
\wabs{\wplagrxu{k}{\hat{x}}{\hat{x}} \wlr{\hat{x} - \hat{x}_k}}.
\end{equation}
\peTheorem{thmLipschitz}
The next lemma shows that the term $\wfc{\tau}{\wvec{\hat{x}}}$ in Theorem \ref{thmLipschitz} is not much
different from $\wfc{\tau}{\wvec{x}}$ when $\Delta$ and $\wnorm{\wvec{z}}_\infty$ are small:
\pbLemma{lemBoundTau}
Under the conditions \pRef{first_cond}--\pRef{last_cond}, consider
$\Delta$ in \pRef{big_delta} and $\wvec{z} = \wfc{\wvec{z}}{\hat{\wvec{x}}, \wfc{\lambda}{\wvec{x}}}$
in \pRef{def_zk}. If $\Delta < 1$ and
$\wnorm{\wvec{z}}_\infty < 1$ then
\[
\wfc{\tau}{\wvec{\hat{x}}} \leq \frac{\wfc{\tau}{\wvec{x}}}{\wlr{1 - \Delta}\wlr{1 - \wnorm{\wvec{z}}_\infty}}.
\]
\peLemma{lemBoundTau}

The bounds above and the results in \cite{SALZER} lead to the following theorem:
\pbTheorem{thmResumeLip}
Under the conditions \pRef{first_cond}--\pRef{last_cond}, let $f: \wrone{} \rightarrow \wrone{}$ be a function with Lipschitz constant $L$
and consider $\wvec{y} = \wfc{f}{\hat{\wvec{x}}}$. In Salzer's case, for all $\hat{x} \in [-1,1]$, the second barycentric formula $q$ in \pRef{second_formula}
satisfies
\begin{equation}
\label{bound_resume_lip}
\wabs{\wfc{q}{\hat{x};\wvec{\hat{x}}^c,\wvec{y}, \wfc{\lambda}{\wvec{x}^c}} - \wfc{P}{\hat{x};\wvec{\hat{x}}^c,\wvec{y}}}
 \leq 1.3 \times 10^{-15} \ L \ n \, \wlr{2.9 + \log n}^2.
\end{equation}
Moreover, for every polynomial $Q$ with degree at most $n$ and $M$ in \pRef{def_q}, we have
\begin{equation}
\label{bound_resume_lip_b}
\wabs{\wfc{q}{\hat{x};\wvec{\hat{x}}^c,\wvec{y}, \wfc{\lambda}{\wvec{x}^c}} - \wfc{P}{\hat{x};\wvec{x}^c,\wvec{y}}}
 \leq 1.3 \times 10^{-15} \ L \ n \, \wlr{2.9 + \log n}^2 + M \wlr{1.4 \log n + 2}.
\end{equation}
\peTheorem{thmResumeLip}
On the other hand,  Theorem \ref{thmLipschitz},
Lemma \ref{lemSalzersConstants}, and equation \pRef{lip_bound} lead to the following:
\pbTheorem{thmEmpirical}
Under the conditions \pRef{first_cond}--\pRef{last_cond}, let $f: \wrone{} \rightarrow \wrone{}$ be a function with Lipschitz constant $L$
and consider $\wvec{y} = \wfc{f}{\hat{\wvec{x}}}$.
In Salzer's case, with $n \leq 10^6$, and the set $C$ formed by the rounded Chebyshev points corresponding to $2n$,
for all $\hat{x} \in [-1,1]$ the second barycentric formula $q$ in \pRef{second_formula}
satisfies
\begin{equation}
\label{bound_empirical}
\wabs{\wfc{q}{\hat{x};\wvec{\hat{x}}^c,\wvec{y}, \wfc{\lambda}{\wvec{x}^c}} - \wfc{P}{\hat{x};\wvec{\hat{x}}^c,\wvec{y}}}
\leq
0.8 \, L \ \wfc{b_n}{C,\wvec{z}} \, \wlr{2.3 + \log n}.
\end{equation}
Moreover, for every polynomial $Q$ with degree at most $n$ and $M$ in \pRef{def_q}, we have
\[
\wabs{\wfc{q}{\hat{x};\hat{x}^c,\wvec{y}, \wfc{\lambda}{\wvec{x}^c}} - \wfc{P}{\hat{x};\wvec{x}^c,\wvec{y}}}
\leq
0.8 \, L \ \wlr{\wfc{b_n}{C,\wvec{z}} + 3.7 \times 10^{-16}}  \, \wlr{2.3 + \log n}
\]
\begin{equation}
\label{bound_empirical_b}
+ M \wlr{1.4 \log n + 2}.
\end{equation}
\peTheorem{thmEmpirical}
In order to use Theorem \ref{thmEmpirical}, we may need to rely on empirical data regarding $b_n$, like the data
presented in Table \ref{table_bn}. However, Theorem \ref{thmEmpirical} yields a sharper upper bound than
Theorem \ref{thmResumeLip} for large $n$. For instance, we have $b_n \leq 2 \times 10^{-15}$
in all entries in Table \ref{table_bn}, and no experiment we performed resulted
in $b_n$ greater than $3 \times 10^{-15}$.
Therefore, our empirical data in combination with Theorem \ref{thmEmpirical} suggest
that, in Salzer's case, the forward error in Step II for the second formula will be at most
$2.2 \times 10^{-15} L \wlr{2.3 + \log n}$ for $n$ up to $10^6$, and this bound is much smaller than
the one provided by Theorem \ref{thmResumeLip}.

Theorems \ref{thmResumeLip} and \ref{thmEmpirical} bound the forward error in Step II for
the second formula, and the next theorem presents a bound on the forward error on Step III. By combining
these bounds with the bounds above we obtain an overall bound for the numerical forward errors for
the second formula.

\pbTheorem{thmStepIIISec}
Under the conditions \pRef{first_cond}--\pRef{last_cond},
let $\epsilon$ be the machine precision and assume that
$\wlr{2 n + 6} 
\epsilon < 1$ and 
$\wnorm{\wfc{\bm{\zeta}}{\wfc{\lambda}{\hat{\wvec{x}}},\hat{\wvec{w}}}}_\infty \wlr{1 + \Lambda_{\hat{x}^-,\hat{x}^+,\hat{\wvec{x}}}} < 1$
and define
\[
\Lambda :=
\frac{\Lambda_{\hat{x}^-,\hat{x}^+,\hat{\wvec{x}}}}{1 - \wnorm{\wfc{\bm{\zeta}}{\wfc{\lambda}{\hat{\wvec{x}}},\hat{\wvec{w}}}}_\infty \wlr{1 + \Lambda_{\hat{x}^-,\hat{x}^+,\hat{\wvec{x}}}} }.
\]
If $\wlr{n + 2} \wlr{2 + \Lambda} \epsilon < 1$
then for every $\hat{x} \in [\hat{x}^-, \hat{x}^+]$ and $\wvec{y} \in \wrn{n+1}$
the computed value $\wrounde{\wfc{q}{\hat{x}, \hat{\wvec{x}}, \wvec{y}, \hat{\wvec{w}}}}$ satisfies
\begin{equation}
\label{bound_stepIII_a}
\wabs{\wrounde{\wfc{q}{\hat{x}, \hat{\wvec{x}}, \wvec{y}, \hat{\wvec{w}}}} -
\wfc{q}{\hat{x}, \hat{\wvec{x}}, \wvec{y}, \hat{\wvec{\wvec{w}}}}}
\leq \Lambda \frac{\theta + 2 n + 6}{1 - \wlr{2 n + 6} \epsilon} \wnorm{\wvec{y}}_\infty  \epsilon,
\end{equation}
for
\[
\theta := \frac{\wlr{1 + \Lambda} \wlr{n + 2}}{1 - \wlr{n + 2} \wlr{2 + \Lambda} \epsilon }.
\]
In Salzer's case with $\epsilon \leq 2.3 \times 10^{-16}$, the
forward error in Step III is bounded by
\[
\wabs{\wrounde{\wfc{q}{\hat{x}, \hat{\wvec{x}}^c, \wvec{y}, \wfc{\lambda}{\wvec{x}^c}}} -
\wfc{q}{\hat{x}, \hat{\wvec{x}}, \wvec{y}, \wfc{\lambda}{\wvec{x}^c}}}
\]
\begin{equation}
\label{bound_stepIII_b}
\leq 0.26 \wlr{2 n + 6} \wlr{1.6 + \log n} \wlr{5.8 + \log n} \wnorm{\wvec{y}}_\infty \epsilon \approx 0.5 \, n \, \wfc{\log^2}{n} \wnorm{\wvec{y}}_\infty \epsilon.
\end{equation}
\peTheorem{thmStepIIISec}
Finally, we note that \cite{Masc} also presents bounds for Step III applicable in Salzer's case.
On the one hand, some bounds in that article involve the Lipschitz constant, on the
other hand, they do not have $\wfc{O}{\log n}$ factors.

\section{Improving the accuracy of the first formula in Salzer's case}
\label{section_new_stuff}
The results in the previous sections show that the accuracy of both barycentric
formulae in the Step II in Figure \ref{figure_steps} is affected by the relative
errors in the length of the intervals $[x_j,x_k]$ (which we measure by
$\delta$ and $\Delta$) and lead to the relative errors $z_k$ in the weights.
Moreover, Table \ref{table_rho} and Figure \ref{figure_least_squares} show that
the first formula is quite sensitive to errors in this second step.

In Salzer's case, when $n$ is large, the $\delta_{k(k+1)}$ for $k$ near $0$ are much larger than
the $\delta_{k(k+1)}$ for $k$ near $n/2$, because the intervals $[x_k,x_{k+1}]$ for
$k$ small have length of order $1/n^2$ whereas the intervals
$[x_k,x_{k+1}]$ for $k$ near $n/2$ have length of order $1/n$.
As a result, the $z_k$s tend also to be larger for $x_k$ near $\pm 1$.
Similar discrepancies in the sizes of the $\delta_{jk}$ and the $z_k$ will happen whenever we
round the nodes as usual and they cluster around a point $c \in [x^-,x^+]$.
This section proposes an efficient way to reduce the
largest $z_k$, by improving the accuracy of the nodes near
the points at which they cluster, and presents experiments
showing that this procedure is effective for the first formula
in Salzer's case. In principle, our strategy would work for any family
of nodes, and any accumulation point $c$, but it is particularly
appropriate for the Chebyshev points,
because in this case the trigonometric identity
$1 - \wfc{\cos}{k\pi/n} = 2 \wfc{\sin^2}{k \pi/\wlr{2n}}$
allows us to evaluate the difference $x_k - x_0$ accurately in
double precision, and the difference $x - x_0 = 1 + x$ is computed
exactly for a floating point number $x \in [-1,-1/2)$. In other situations,
for a general cluster point $c$, the implementation
of our strategy may be more difficult.

 This section is divided into three subsections.
 In the first subsection, we describe a new finite precision representation of the nodes and we explain why it improves the accuracy of the first formula.
We then present experimental results showing that our finite precision representation of the nodes is practical in terms of performance and accuracy.
We conclude with some remarks about the new node representation.

\subsection{A new finite precision representation for the nodes}
We now describe a finite precision representation of the nodes that
lead to smaller $\delta_{jk}$ and $z_k$ without the use of quadruple precision.
We partition the interval $[\hat{x}^-,\hat{x}^+]$ into sub-intervals which we will refer to as {\it bins}.
The $l$-th bin has a base point $b_l$ and the nodes $x_k$ in this
bin are represented using $r_k := x_k - b_l$, so that
 $x_k = b_l + r_k$.
We store rounded versions $\hat{r}_k$ of $r_k$ instead of $\hat{x}_k$.
This idea is similar to Dekker's \cite{DEKKER} but our approach is more economical since  we only store one $b_l$ per bin.
For example, we could
use the bins $[-1,-0.5)$, $[-0.5,0.5]$ and $(0.5,1]$
as  illustrated in  Figure \ref{figure_three_bins}:

\begin{center}
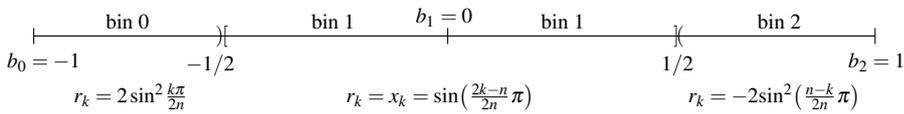
\begin{figure}[!h]
\begin{picture}(0,35)(0,5)
\put(10,30){\line(1,0){315}}
\put(10,27){\line(0,1){6}}
\put(0,18){$b_0 = -1$}
\put(25,5){$r_k = 2 \sin^2 \! \frac{k \pi}{2n}$}

\put(78,28){$)$}
\put(80.5,28){$[$}
\put(67,17){$-1/2$}
\put(114,33){bin 1}
\put(37,33){bin 0}

\put(165,27){\line(0,1){6}}
\put(153,35){$b_1 = 0$}
\put(127,5){$r_k = x_k = \wfc{\sin}{\frac{2 k - n}{2n}  \pi}$}

\put(250.5,28){$($}
\put(249.5,28){$]$}
\put(245,17){$1/2$}
\put(281,33){bin 2}
\put(200,33){bin 1}

\put(325,27){\line(0,1){6}}
\put(315,18){$b_2 = 1$}
\put(255,5){$r_k = - 2 \wfc{\sin^2}{\frac{n - k}{2n}  \pi}$}
\end{picture}
\caption{Partitioning $[-1,1]$ into bins with base points $b_0 = -1$, $b_1 = 0$ and $b_2 = 1$.}
\label{figure_three_bins}
\end{figure}
\end{center}

The $r_k$ in bins $0$ and $2$ can be
computed with much smaller absolute errors than the corresponding $\hat{x}_k$.
The relative errors in $\hat{r}_k$ and $\hat{x}_k$ have the same order of magnitude.
However, $\hat{r}_k$ is smaller than $\hat{x}_k$ and the absolute error in
$\hat{r}_k$ for $x_k$ near the border is orders of magnitude smaller than
the error in $\hat{x}_k$ for large $n$.

The key point in the representation $x_k = b_l + r_k$ is the possibility
of computing the differences $x - x_k$ without evaluating $x_k$
explicitly.
Instead, given $x \in [\hat{x}^-,\hat{x}^+]$,
\begin{enumerate}
\item[(a)] we find $l$ such that $x$ is in the $l$-th bin and define $r_x := x - b_l$.
\item[(b)] Given a node $x_k$ in the $m$-th bin, we evaluate $x - x_k$ as $\wlr{b_l - b_m} + (r_x - r_k)$.
\end{enumerate}

The steps (a) and (b) are accurate because:
\begin{enumerate}
\item[(i)] We choose the bins and their bases so that we can apply the
following lemma
and conclude that there is no rounding errors in the evaluation of $r_x$ in step (a).
\pbLemma{lemSterbenz}(Sterbenz's Lemma)
If subtraction is performed with a guard digit and $x/2 \leq y \leq 2 x$ then $x - y$ is computed exactly.
\peLemma{lemSterbenz}
\item[(ii)] We choose base points $b_l$ such that the differences $b_l - b_m$ are computed exactly.
\item[(iii)] When $x$ and $x_k$ are very close, the difference $b_l - b_m$ is
 small and the difference in step (b) is computed with high relative accuracy.
\end{enumerate}

The statements above can be formalized as follows:
\pbTheorem{thmGoodRik}
In Salzer's case, with $x_k = b_l + r_k$ for $r_k$ as in Figure \ref{figure_three_bins},
 consider the nodes $\tilde{x}_k = b_l + \hat{r}_k$. If $\hat{r}_k = r_k \wlr{1 + \theta_k}$, with $\wabs{\theta_k} \leq 4.6 \times 10^{-16}$,
 then the numbers $\wfc{\delta_{jk}}{\tilde{\wvec{x}},\wvec{x}}$ in \pRef{def_delta} satisfy
\begin{equation}
\label{good_rik}
\sum_{j \neq k} \wabs{\wfc{\delta_{jk}}{\tilde{\wvec{x}},\wvec{x}}} \leq
\wnorm{\bm{\theta}}_\infty \, \wlr{3.2 + 2.3 \, n + 4.3 \, n \,\wfc{\log}{n+ 1}}.
\end{equation}
\peTheorem{thmGoodRik}
\noindent
The bound \pRef{good_rik} is of order $\wfc{O}{\epsilon n \log n}$ whereas the analogous
bound for nodes rounded as usual $\wvec{z}$ is
of order $\wfc{\Theta}{\epsilon n^2}$. This explains why bins
improve the accuracy of the first formula in Salzer's case. The bound \pRef{good_rik}
and Lemma \ref{lemBoundZk} lead to smaller
upper bounds for the forward errors for the second formula in Theorem \ref{thmMain}
and to smaller bounds on the backward errors for the second formula presented in
\cite{MascCam}.

\subsection{Nodes in bins are practical}
\label{subsection_experiments}
 In this section we compare common implementations for both barycentric formulae with implementations based on bins (as described in Appendix \ref{section_experiments}).
Our results are presented in tables
\ref{table_stepII_error}, \ref{table_overall_error} and \ref{table_times}.
Table \ref{table_stepII_error} shows that by using bins we reduce
the errors introduced by Step II for the first formula.  However, even with this reduction,
the errors introduced by Step II for the first formula are larger
than those for the second formula.
 Table \ref{table_stepII_error} also shows that using bins did not improve the second formula,
 because the Step II errors with usual nodes for this formula are already small.

\begin{table}[!h]
\caption{Maximum error in Step II for $\wfc{f}{x} = \wfc{\cos}{100x}$ and $x \in X_{-1,n}$.}
\centering
\begin{tabular}{c|c|ccc|c|ccc}
\hline\\[-0.27cm]
      & \multicolumn{4}{c|}{First Formula}  &  \multicolumn{4}{c}{Second Formula} \\[-0.01cm]
\cline{2-9}\\[-0.25cm]
$n + 1$ & as usual & 3 bins & 39 bins & 79 bins & as usual  & 3 bins & 39 bins & 79 bins  \\[-0.01cm]
\hline\\[-0.2cm]
$ 10^3$ & 8.9e-12 & 	5.0e-15 & 2.1e-15 & 2.0e-15 & 4.4e-16 & 9.7e-17 & 6.8e-17 & 6.8e-17\\
$ 10^4$ & 4.4e-10 & 	7.8e-15 & 9.0e-15 & 2.7e-15 & 9.0e-17 & 9.0e-17 & 9.0e-17 & 9.0e-17\\
$ 10^5$ & 7.1e-08 & 	8.5e-15 & 7.8e-15 & 3.1e-15 & 8.0e-17 & 7.7e-17 & 7.7e-17 & 7.7e-17\\
$ 10^6$ & 2.1e-06 & 	1.1e-14 & 1.6e-14 & 1.2e-14 & 7.2e-17 & 7.7e-17 & 7.7e-17 & 7.7e-17
\end{tabular}
\label{table_stepII_error}
\end{table}

 Table \ref{table_overall_error} considers the error in the three steps in Figure 1 and shows that, indeed,
bins make the first formula as accurate as the second overall.
Table \ref{table_times} shows that, in these experiments, nodes in bins lead to the same performance
as the usual ones.
These results can be explained by the need for only  one extra sum per node in bin and that the dominant factors in performance are access to memory and divisions.
 In order to evaluate  $\Delta x_k := x - x_k$ for all $x_k$ in the $m$-th bin when $x$ is in the
 $l$-th bin, we compute $\Delta b = b_l - b_m$
 only once and then compute $\Delta \! x_k = \Delta b + (r_x - r_k)$.
 The cost incurred in evaluating this expression is less than
 twice the cost incurred in evaluating $\Delta \! x_k$ as usual because $\Delta b$ stays in the cache.

\begin{table}[!h]
\caption{Maximum overall error for $\wfc{f}{x} = \wfc{\cos}{100x}$ and $x \in X_{-1,n}$.}
\centering
\begin{tabular}{c|c|ccc|c|ccc}
\hline\\[-0.27cm]
      & \multicolumn{4}{c|}{First Formula}  &  \multicolumn{4}{c}{Second Formula} \\[-0.01cm]
\cline{2-9}\\[-0.22cm]
$n + 1$    & as usual & 3 bins & 39 bins & 79 bins & as usual  & 3 bins & 39 bins & 79 bins  \\[-0.01cm]
\hline\\[-0.2cm]
$ 10^3$ & 8.9e-12 & 1.2e-14 & 1.2e-14 & 1.2e-14 & 1.1e-14 & 1.1e-14 & 1.0e-14 & 1.0e-14 \\
$ 10^4$ & 4.4e-10 &	3.5e-14 & 3.3e-14 & 3.1e-14 & 3.0e-14 & 3.0e-14 & 3.0e-14 & 3.0e-14 \\
$ 10^5$ & 7.1e-08 & 8.9e-14 & 9.6e-14 & 9.6e-14 & 8.9e-14 & 8.7e-14 & 8.7e-14 & 8.7e-14 \\
$ 10^6$ & 2.1e-06 & 2.4e-13 & 2.5e-13 & 2.4e-13 & 1.7e-13 & 1.6e-13 & 1.6e-13 & 1.6e-13
\end{tabular}
\label{table_overall_error}
\end{table}

\begin{table}[!h]
\caption{Normalized times (usual first formula = 1).}
\centering
\begin{tabular}{c|c|ccc||c|ccc}
\hline\\[-0.27cm]
      & \multicolumn{4}{c|}{First Formula}  &  \multicolumn{4}{c}{Second Formula} \\[-0.01cm]
\cline{2-9}\\[-0.22cm]
$n + 1$ & as usual & 3 bins & 39 bins & 79 bins & as usual & 3 bins & 39 bins & 79 bins  \\[-0.01cm]
\hline\\[-0.25cm]
$ 10^3$ & $1.00$      &   1.01 & 1.05 &  1.07   & 1.40     & 1.40  & 1.40     &  1.40    \\
$ 10^4$ & $1.00$      &   1.00 & 1.01 &  1.01   & 1.40     & 1.41  & 1.41     &  1.40    \\
$ 10^5$ & $1.00$      &   0.99 & 0.99 &  0.99   & 1.40     & 1.39  & 1.39     &  1.39    \\
$ 10^6$ & $1.00$      &   1.00 & 1.00 &  1.00   & 1.38     & 1.38  & 1.37     &  1.38    \\
\end{tabular}
\label{table_times}
\end{table}

Tables \ref{table_overall_error} and \ref{table_times}
indicate that, in these experiments, the first formula with bins is
competitive in terms of performance and accuracy, whereas the accuracy of the
usual first formula is unacceptable.
 We doubt that the accuracy  will be affected significantly with the use of other compilers, machines and programming languages.
 The performance, on the other hand, depends
on the machine, the compiler and the language. For instance,
 \verb Fortran  sometimes
leads
to faster code than \verb C++,  and by using
a different language we could obtain a different relation among the performance of the
several alternatives described in the previous tables.

\subsection{Final remarks}
The approximate nodes $\tilde{x}_k$ in Theorem \ref{thmGoodRik} may not be floating
point numbers, in the same way that Dekker's numbers usually are not floating
point numbers.
This is not a problem if we already have the $y_k$. In this case, we do
not need the $\tilde{x}_k$, because both formulae can be expressed
in terms of the $y_k$ and the differences $x - \tilde{x}_k$,
and such differences can be computed without the explicit value of $\tilde{x}_k$.
However, any $\tilde{x}_k$ that is not a floating point number
 complicates the evaluation of $y_k = \wfc{f}{\tilde{x}_k}$.
We could handle this problem in three ways:
\begin{itemize}
\item Evaluate $\wfc{f}{\hat{x}_k}$, where $\hat{x}_k$ is the
floating point number closest to $\tilde{x}_k$.
\item Estimate $\wfc{f}{\tilde{x}_k}$ as
$\wfc{f}{r_k} + \wdfc{f}{\hat{x}_k} \wlr{\tilde{x}_k - \hat{x}_k}$ when $f$ is differentiable.
\item Evaluate $\wfc{f}{\tilde{x}_k}$, on-line or off-line, using higher precision arithmetic.
\end{itemize}
In any case, the readers will need to take
these considerations into account should
they  decide to apply the ideas presented in this section.

There are many choices for the bins and their base points
but there isn't a single choice that is optimal
for all compilers, processors and instruction sets.
For instance, the advances in hardware
may lead to efficient combinations of integer and floating point arithmetic.
 We experimented with $\hat{x}_k$
represented as 64 bit integers and the resulting
code was  50\% slower than the one using only floating point arithmetic.
However, with integer arithmetic we reduced the Step II errors by a factor
of $10^3$, due to a better use of the $11$ bits that represent the exponent in IEEE754 double
precision.

The errors in Step II for the implementations
described in the previous subsection are much smaller than the Step III errors.
As a result, the overall error is determined by Step III and the gain
in accuracy in Step II due to the use of integer arithmetic was irrelevant.
Therefore, our mixing of floating point and integer arithmetic is not competitive at
this time.
However, parallel usage of the hardware dedicated to floating point and integer arithmetic could change this.

Finally, in extreme situations, bins can also improve the accuracy of
the second formula. In our experiments with $\wfc{f}{x} = \wfc{\cos}{10^5 x}$ and $10^6$
nodes, the maximum error with the usual second formula was of the order $10^{-12}$ and the
error with the second formula with nodes in bins was of the order $10^{-13}$. Moreover,
the results in section 3 of \cite{MascCam} indicate that for some functions with large Lipschitz
constants the backward error for the second formula could also be reduced by using bins.

\appendix
\section{Proofs}
\label{section_proofs}
This appendix proves the results stated in the previous sections. We state four
more lemmas, after that we prove all lemmas in the order in which they were
stated, we then prove the theorems, also in the order in which they were stated.


\pbLemma{lemBoundProd}
Given a vector $\wvec{v} = \wlr{v_0,v_1 \dots v_n}^t \in \wrn{n+1}$, define
$s := \sum_{k=0}^n v_k$,  $\ a := \sum_{k = 0}^n \wabs{v_k}$,
\[
s_- :=  \wfc{s_-}{\wvec{v}} := - \sum_{k = 0}^n \min \wset{v_k,0}.
\hspace{1cm} \wrm{and} \hspace{1cm}
s_+ := \wfc{s_+}{\wvec{v}} := \sum_{k = 0}^n \max \wset{v_k,0}.
\]
If $s_- < 1$ then
\begin{equation}
\label{bound_inv_prod}
-s \leq \prod_{k=0}^n \frac{1}{1 + v_k} - 1 \leq \frac{s_-}{\wlr{1 - s_-}\wlr{1 + s_+}} - \frac{s_+}{1 + s_+} \leq \frac{s_-}{1 - s_-}
\end{equation}
and
\begin{equation}
\label{lower_bound_prod}
\prod_{k =0}^n \wlr{1 + v_k} - 1 - s  \geq - s_- s_+ \geq - \frac{a^2}{4}.
\end{equation}
Moreover, if $s_- < 1$ and $s < 1$ then
\begin{equation}
\label{upper_bound_prod}
\prod_{k =0}^n \wlr{1 + v_k} - 1 - s  \leq \frac{s^2}{1 - s}.
\end{equation}
Finally, if $a < 1$, then, for $0 \leq m \leq n$, the product
\[
P := \wlr{\prod_{k = 0}^m \frac{1}{1 + v_k}} \wlr{\prod_{k= m+1}^n \wlr{1 + v_k}}
\]
satisfies
\begin{equation}
\label{bound_prod_abs}
1 - a \leq P \leq \frac{1}{1- a} \hspace{1cm} \wrm{and} \hspace{1cm} \wabs{P - 1} \leq \frac{a}{1 - a}.
\end{equation}
\peLemma{lemBoundProd}

%

\pbLemma{lemSums}
For $1 \leq k \leq n/2$, the Chebyshev points of the second kind satisfy
\begin{equation}
\label{good_rk_b}
\sum_{j = 0}^{k-1} \frac{1}{x_k - x_j}
\leq \frac{1}{\wfc{2 \sin}{\frac{k \pi}{n}} \wfc{\sin}{\frac{k \pi}{2n}}} \wlr{1 + k \, \wfc{\log}{4k - 1}}
\leq \frac{n^2}{4 \sqrt{2} \ k^2} \wlr{1 + k \, \wfc{\log}{4k - 1}}
\end{equation}
and
\begin{equation}
\label{good_rk_c}
\sum_{j = k+1}^{n} \frac{1}{x_j - x_k} \leq
\frac{3n^2}{4 \sqrt{2}\, k}  \wfc{\log}{4 k + 1}.
\end{equation}
Since \pRef{good_rk_b} and \pRef{good_rk_c} decrease with $k$, for $1 \leq k \leq n/2$
\[
\sum_{j \neq k} \frac{1}{\wabs{x_k - x_j}} \leq
\frac{n^2}{4 \sqrt{2}} \wlr{ 1 + \log \!3 + 3 \log \! 5} < 1.2246\, n^2
\hspace{0.7cm} \wrm{and} \hspace{0.7cm}
\sum_{j=1}^n \frac{1}{x_j -x_0} \leq \frac{\pi^2 n^2}{12} < 0.82247 \, n^2.
\]
In particular, for all $n$,
\pbHypot{goodRkBC}
\sum_{j \neq k} \frac{1}{\wabs{x_k - x_j}} \leq 1.2246\, n^2.
\peHypot{goodRkBC}
Moreover, if $n \geq 10$ then
\pbDef{defMu}
\mu := \frac{1}{\min_{0 \leq k < n} \wlr{x_{k+1}- x_k}} = \frac{1}{\min_{j \neq k} \wabs{x_j - x_k}}
\peDef{defMu}
satisfies
\pbTClaim{boundMu}
\mu \leq 0.20432 \, n^2.
\peTClaim{boundMu}
\peLemma{lemSums}

\pbLemma{lemDelta}
If $x^- = x_0$, $x^+ = x_n$ and $\mu$ in \pRef{defMu}
is such that $2 \mu \wsupdx < 1$,
then $\delta_{jk}$, $\Delta_k$ in
\pRef{def_delta}--\pRef{big_delta} and $a_k := \sum_{j \neq k} \delta_{jk}$ satisfy
\begin{equation}
\label{bound_delta}
\wabs{\delta_{jk}} \leq \frac{\kappa}{\wabs{x_j - x_k}} \ \ \wrm{for} \ j \neq k,
\hspace{2cm}
a_k \leq  \kappa  \sum_{j \neq k} \frac{1}{\wabs{x_j - x_k}}
\end{equation}
and, for $0 \leq k < n$,
\begin{equation}
\label{bound_big_delta}
\Delta_k \leq \kappa \wlr{ \sum_{j = 0}^{k-1} \frac{1}{x_k - x_j}
+ \frac{2}{x_{k+1} - x_k} + \sum_{j = k + 2}^n \frac{1}{x_j - x_{k+1}}},
\end{equation}
for
\pbDef{defKappa}
\kappa := \frac{2 \wsupdx}{1 -  2  \mu \wsupdx }
\peDef{defKappa}
\peLemma{lemDelta}

\pbLemma{lemLipschitz}
Consider nodes $\wvec{x} = \wset{x_0,\dots,n}$ and and interval $[x^-,x^+]$. Given $0 \leq k \leq n$,
define $\wvec{x}' := \wset{x_0,\dots,n} \setminus \wset{x_k}$. If $x \in [x^-,x^+]$ then
\begin{equation}
\label{bound_lipschitz}
\wabs{\wfc{P}{x;\wvec{x},\wvec{y}} - y_k} \leq \Lambda_{x^-,x^+,\wvec{x}'}
\wabs{x - x_j} \max_{j \neq k} \wabs{\frac{y_j - y_k}{x_j - x_k}}
\end{equation}
\peLemma{lemLipschitz}

\subsection{Proofs of the lemmas}

\pbProof{Lemma}{lemBoundZk}
According to definition \pRef{def_lambda}, we have
\[\wfcc{\lambda_k}{\wvec{x}} = \frac{1}{\prod_{j \neq k} \wlr{x_k - x_j}} \hspace{1cm} \wrm{and} \hspace{1cm}
\wfcc{\lambda_k}{\wvec{\hat{x}}} = \frac{1}{\prod_{j \neq k} \wlr{\hat{x}_k - \hat{x}_j}}.
\]
Using definitions \pRef{def_delta} and \pRef{def_zk} we obtain
\[
z_k = \frac{\wfcc{\lambda_k}{\wvec{x}}}{\wfc{\lambda_k}{\wvec{\hat{x}}}} - 1
= \wlr{\prod_{j \neq k} \frac{\hat{x}_k - \hat{x}_j}{x_k - x_j}} \ - 1 =
 \wlr{\prod_{j \neq k} \frac{1}{\frac{x_k - x_j}{\hat{x}_k - \hat{x}_j}}} \ - 1 =
\wlr{\prod_{j \neq k}\frac{1}{1 + \delta_{kj}}} - 1,
\]
and the bounds on $z_k$ in Lemma \ref{lemBoundZk} follow from Lemma \ref{lemBoundProd}.
Similarly,
\[
\zeta_k = \frac{\wfc{\lambda_k}{\wvec{\hat{x}}}}{\wfcc{\lambda_k}{\wvec{x}}} - 1
= \wlr{\prod_{j \neq k} \frac{x_k - x_j}{\hat{x}_k - \hat{x}_j}} \ - 1 =
\wlr{\prod_{j \neq k}\wlr{1 + \delta_{kj}}} - 1,
\]
and the bound in $\wabs{\zeta_k}$ also follows from Lemma \ref{lemBoundProd}
\peProof{Lemma}{lemBoundZk}


\pbProof{Lemma}{lemBackwardProduct}
Consider a set of indices $K \subset \wset{0,1,\dots, n}$.
If $\hat{x} = \hat{x}_k$ for $k \in K$ then we can take $x = x_k$.
If $\hat{x} = \hat{x}_k$ for $k \not \in K$ then by taking $x = x_k$ and using the
definition \pRef{def_delta} we obtain
\[
\wabs{\frac{x - x_j}{\hat{x} - \hat{x}_j} - 1} = \wabs{\delta_{kj}}
\]
for all $j \in K$. On the other hand, if $\hat{x} \in [\hat{x}^-,\hat{x}^+] \setminus \wset{\hat{x}_0,\dots,\hat{x}_n}$
then Corollary 4.2 in \cite{MascCam}
\footnote{Please note that there is a typo in \cite{MascCam}: the indexes $k$ and $j$ are exchanged.}
shows that there exists $x \in [x^-,x^+] \setminus \wset{x_0,\dots,x_n}$
which satisfies \pRef{delta_xu} and, for $0 \leq j \leq n$,
\begin{itemize}
\item[(a)] If $\hat{x}^- < \hat{x} < \hat{x}_{k^-}$ then
\[
\wabs{\frac{x - x_j}{\hat{x} - \hat{x}_j} - 1} \leq
\max \wset{ \wabs{\delta^-_{j}}, \ \wabs{\delta_{kj^-}}}.
\]
\item[(b)] If $k^- \leq k <  k^+$ and $\hat{x}_k < \hat{x} < \hat{x}_{k + 1}$ then
\[
\wabs{\frac{x - x_j}{\hat{x} - \hat{x}_j} - 1} \leq
\max \wset{\wabs{\delta_{kj}}, \ \wabs{\delta_{k\wlr{j+1}}}}.
\]
\item[(c)] If $\hat{x}_{ k^+} < \hat{x} < \hat{x}^+$ then
\[
\wabs{\frac{x - x_j}{\hat{x} - \hat{x}_j} - 1} \leq
\max \wset{\wabs{\delta_{k j^+}}, \ \wabs{\delta^+_j}},
\]
\end{itemize}
where $k^-$ and $k^+$ are such that
$\hat{x}_{k^+}$ is the first node largest than $\hat{x}^-$ and
$\hat{x}_{k^+}$ is the last node smaller than $\hat{x}^+$.
Defining $d_j = \wlr{x - x_j}/\wlr{\hat{x} - \hat{x}_j} - 1$,
we can write the product in the left hand side of \pRef{backward_product} as
\[
\prod_{j \in K} \wlr{\hat{x} - \hat{x}_j} =
\wlr{\prod_{j \in K} \wlr{\frac{\hat{x} - \hat{x}_j}{x - x_j}}} \prod_{j \in K} \wlr{x - x_j} =
\wlr{\prod_{j \in K} \frac{1}{1 + d_j}} \prod_{j = 0}^n \wlr{x - x_j} =
\wlr{1 + \beta} \prod_{j = 0}^n \wlr{x - x_j}
\]
for
\[
\beta := \prod_{j \in K} \frac{1}{1 + d_j} - 1.
\]
In case (a) above we have $\sum_{j \in K} \wabs{d_j} \leq \Delta^-$,
in case (b) or when $\hat{x} = \hat{x}_k$ for $k \not \in K$ we have $\sum_{j \in K} \wabs{d_j} \leq \Delta_k$ and in case (c)
we have $\sum_{j \in K} \wabs{d_j} \leq \Delta^+$.
In all cases Lemma \ref{lemBackwardProduct} follows from Lemma \ref{lemBoundProd}
\peProof{Lemma}{lemBackwardProduct}

\pbProof{Lemma}{lemSalzersConstants}
\pbClaim{hsc}
Recall that in \hyperlink{link_salzer}{Salzer's case} we are
are constrained by \pRef{salzersCase}:
\pbTClaim{hsc}
10 \leq n \leq 2 \times 10^{6}
\hspace{1cm} \pStop{and} \hspace{1cm}
\wsupdx \leq 4.6 \times 10^{-16}.
\peEClaim{hsc}
\pbClaim{boundKappa}
Equations \pRef{boundMu} and \pRef{hsc} show that the constant
\lpBind{$\kappa$}{in \pRef{defKappa}}
is bounded by
\pbChain{boundKappaA}
\kappa \leq
\frac{2 \wsupdx}{1 - 2 \wsupdx \times 0.20432 \times n^2}
\leq
\peChain{boundKappaA}
\pbTClaim{boundKappa}
\leq \frac{2 \wsupdx}{1 - 9.2 \times 10^{-16} \times 0.20432 \times 4 \times 10^{12}}
\leq 2.0016 \wsupdx,
\peEClaim{boundKappa}
\hypertarget{at_work}{}
\pbClaim{sharpDelta}
and equations \pRef{bound_delta}, \pRef{hsc} and \pRef{boundKappa} yield
\pbChain{sharpDeltaA}
\delta \leq \kappa \mu \leq 2.0016 \times 0.20432 \wsupdx n^2 \leq
\peChain{sharpDeltaA}
\pbTClaim{sharpDelta}
\leq 0.40897 \wnorm{\wvec{x} - \wvec{\hat{x}}}_\infty n^2 \leq 1.8813 \times 10^{-16} n^2 \leq 7.5252 \times 10^{-4}.
\peEClaim{sharpDelta}
\pbClaim{boundGammaK}
Equations \pRef{goodRkBC}, \pRef{bound_delta}, \pRef{hsc} and \pRef{boundKappa}
lead to the following bound on $a_k = \sum_{j \neq k} \wabs{\delta_{jk}}$:
\pbTClaim{boundGammaK}
a_k \leq \kappa \sum_{j \neq k} \frac{1}{\wabs{x_j - x_k}}
\leq 2.4512 \wnorm{\wvec{x} - \wvec{\hat{x}}}_\infty n^2 \leq  1.1276 \times 10^{-15} n^2 \leq 4.5104 \times 10^{-3}.
\peEClaim{boundGammaK}
\pbClaim{sharpZ}
Equations \pRef{boundZk} and \pRef{boundGammaK} yield
\pbTClaim{sharpZ}
\wabs{z_k} \leq \frac{a_k}{1 - a_k} \leq \frac{2.4512 \wnorm{\wvec{x} - \wvec{\hat{x}}}_\infty n^2}{1 - 4.5104 \times 10^{-3}}
\leq 2.4624 \wnorm{\wvec{x} - \wvec{\hat{x}}}_\infty n^2
\leq 1.1328 \times 10^{-15} n^2 \leq 4.5312 \times 10^{-3}.
\peEClaim{sharpZ}
The same bounds hold for $\zeta_k$, because Lemma \ref{lemBoundZk} also states that
$\wabs{\zeta_k} \leq a_k / \wlr{1 - a_k}$.
Lemmas  \ref{lemDelta} and \ref{lemSums} and equation \pRef{boundKappa} also lead to
\[
\Delta_0 \leq \kappa \wlr{\frac{2}{x_1 - x_0} + \sum_{j = 2}^n \frac{1}{x_j - x_1}} \leq
\kappa \wlr{2 \times  0.20432 \times n^2 + \frac{3n^2}{4 \sqrt{2}} \log 5}
\leq 2.5264 \wsupdx  n^2.
\]
For $1 \leq k \leq n/2$, noting that the second sum in \pRef{bound_big_delta} starts at $k + 2$, we obtain
\[
\Delta_k \leq \frac{2 \kappa}{x_{k+1} - x_k} +  \frac{\kappa n^2}{4 \sqrt{2}} \wlr{\frac{1}{k^2} +
\frac{\wfc{\log}{4 k - 1}}{k} +  3 \frac{\wfc{\log}{4 \wlr{k + 1} + 1}}{k + 1}}.
\]
The right hand side of this equations decreases with $k$. By symmetry, for $0 < k < n$,
\begin{equation}
\label{sharp_big_delta}
\Delta_k \leq 2 \kappa \mu +  \frac{\kappa n^2}{4 \sqrt{2}} \wlr{1 + \log 3 + \frac{3}{2} \log 9}
\leq 2.7267 \wnorm{\wvec{x} - \wvec{\hat{x}}}_\infty n^2 \leq  1.2543 \times 10^{-15} n^2 \leq 5.0172 \times 10^{-3}.
\end{equation}
Let us now prove the bound on $\wnorm{\wvec{z}}_1$. When $1 \leq k \leq n/2$,
Lemmas \ref{lemSums} and \ref{lemDelta} show that
\begin{equation}
\label{bound_gamma_psi}
a_k \leq \frac{\kappa n^2}{4 \sqrt{2}} \psi_k,
\end{equation}
for
\[
\psi_k := \frac{1}{k^2} + \frac{\wfc{\log}{4 k - 1}}{k} +  3 \frac{\wfc{\log}{4 k + 1}}{k}.
\]
When $\wabs{x} < 1$, $\wfc{\log}{1 - x} + \wfc{\log}{1 + x} = \wfc{\log}{1 - x^2} < 0$ and
$\wfc{\log}{1 + x} \leq x$. Therefore, defining
\[
\wfc{f}{x} := \frac{3}{2x^2} + \frac{8 \log 2 + 4 \log x}{x},
\]
we obtain
\begin{equation}
\label{bound_psik}
\wfc{\log}{1 - 1/{4k}} + 3 \wfc{\log}{1 + 1/{4k}} \leq \frac{1}{2k}
\hspace{1.0cm}
\wrm{and}
\hspace{1.0cm}
\psi_k \leq \wfc{f}{k}.
\end{equation}
The function $f$ has derivative
\[
\wdfcc{f}{x} = - \frac{3}{x^3} - \frac{8 \log 2 - 4}{x^2} - 4 \frac{\log x}{x^2},
\]
which is negative for $x \geq 1$. Therefore,
\[
\sum_{k = 2}^{\wfloor{n/2}} \psi_k \leq \int_1^{n/2} \wfc{f}{x} \wdx{x} = \left.
- \frac{3}{2 x} + 8 \log 2 \log x + 2 \log^2 x \right|_{x = 1}^{x = n/2}
\]
\begin{equation}
\label{sum_psi_k}
= - \frac{3}{n} + 8 \log 2 \wlr{\log n - \log 2} + 2 \wlr{ \log n - \log 2}^2 + \frac{3}{2}
\leq 2 \wlr{\log n \wlr{\log n + 2 \log 2} + \frac{3}{4} - 3 \log^2 2}.
\end{equation}
By symmetry, and equations \pRef{bound_gamma_psi} and \pRef{sum_psi_k}, we have
\begin{equation}
\label{symmetry}
\sum_{k = 0}^n a_k \leq 2 \wlr{a_0 + a_1} +
\frac{\kappa n^2}{2 \sqrt{2}} \sum_{k = 2}^{\wfloor{n/2}} \psi_k
\leq A \wnorm{\wvec{x} - \wvec{\hat{x}}}_\infty n^2,
\end{equation}
for
\[
A := \frac{2 \wlr{a_0 + a_1}}{\wnorm{\wvec{x} - \wvec{\hat{x}}}_\infty n^2} +
1.0008 \sqrt{2} \wlr{\log n \wlr{\log n + 2 \log 2} + \frac{3}{4} - 3 \log^2 2}.
\]
Using bound \pRef{boundGammaK}, the identity $\wlr{3 + \log n}^2 = \log^2 n + 2 \times 3  \times \log n + 3^2$,
and the hypothesis $n \geq 10$ we deduce that
\[
A \leq 1.0008 \sqrt{2} \wlr{3 + \log n}^2 -
 1.0008 \sqrt{2} \wlr{\frac{33}{4} + 3 \log^2 2 + \wlr{6 - 2 \log 2} \times \log 10}
 +\frac{2 \wlr{a_0 + a_1}}{\wnorm{\wvec{x} - \wvec{\hat{x}}}_\infty n^2}
\]
\[
< 1.0008 \sqrt{2} \wlr{3 + \log n}^2 < 1.415345 \wlr{3 + \log n}^2.
\]
Moreover, \pRef{boundZk} and \pRef{boundGammaK}  show that
$\wabs{z_k} \leq a_k/ \wlr{1 - a_k} \leq 1.0046 \, a_k$.
This inequality, the bound on $A$ above, and \pRef{symmetry} yield
\begin{equation}
\label{sharpZ1}
\wnorm{\wvec{z}}_1 \leq 1.42186 \wnorm{\wvec{x} - \wvec{\hat{x}}}_\infty \wlr{3 + \log n}^2 n^2
\leq 6.5406 \times 10^{-16} \wlr{3 + \log n}^2 n^2 \leq 0.80202,
\end{equation}
and this verifies the last line it table \ref{table_salzers_bounds}
\peProof{Lemma}{lemSalzersConstants}


\pbProof{Lemma}{lemLebesguePerturbation}
In the notation of \cite{MascCam}, we write the $\Lambda_{x^-,x^+,\wvec{x}}$ in \pRef{def_lebesgue} as
$\Lambda_{x^-,x^+,\wvec{x}, \wfc{\lambda}{\wvec{x}}}$. Moreover, if we take
$\wvec{w} = \wfc{\lambda}{\wvec{x}}$ and $\hat{\wvec{w}} = \wfc{\lambda}{\wvec{\hat{x}}}$
then \pRef{def_zk} implies that
$\bm{\zeta} = \wfc{\bm{\zeta}}{\wvec{w},\hat{\wvec{w}}} = \wfc{\wvec{z}}{\hat{\wvec{x}},\wfc{\lambda}{\wvec{x}}}$
and the hypothesis \pRef{hypo_lebesgue_perturbation} shows that $d = \delta$ and
$\bm{\zeta}$ satisfy the hypothesis in Theorem 4.3 in \cite{MascCam}. Therefore,
\pRef{bound_lebesgue} follows from Corollary 3 and Theorem 3 in that article.
Theorem 1 in \cite{GUNTTNER} combined with the results in \cite{MCCABE} yield that,
for the Euler-Mascheroni constant $\gamma < 0.577215665$ and $n \geq 10$,
\begin{equation}
\label{mascheroni}
\Lambda_{-1,1,\wvec{x}^c} \leq \frac{2}{\pi} \wlr{\log n + \gamma + \log \frac{8}{\pi} + \frac{\pi^2}{144 n^2}}
\leq 0.63662 \wlr{\log n + 1.5127}.
\end{equation}
Moreover, evaluating the term after the first $\leq$ in the expression above for $n = 2 \times 10^6$ we obtain
\begin{equation}
\label{leb_salzer_max}
\Lambda_{-1,1,\wvec{x}^c} \leq 10.1991.
\end{equation}
The bounds on $\delta$ and $\wvec{z}$ in Lemma \ref{lemSalzersConstants} and
\pRef{bound_lebesgue} lead to
\begin{equation}
\label{leb_salzer_aux}
\Lambda_{-1,1,\hat{\wvec{x}}^c} \leq 1.0629 \Lambda_{-1,1,\wvec{x}^c}.
\end{equation}
This bound and \pRef{mascheroni} lead to $\Lambda_{-1,1,\hat{\wvec{x}}^c} \leq 0.67667 \log n + 1.0236$,
and \pRef{leb_salzer_max} and \pRef{leb_salzer_aux} show that $\Lambda_{-1,1,\hat{\wvec{x}}^c} \leq 10.841$
\peProof{Lemma}{lemLebesguePerturbation}

\pbProof{Lemma}{lemWarmUp}
Lemma \ref{lemBackwardProduct} yields $\beta$ satisfying
\pRef{bound_warm_up} and equation \pRef{cheby_prod} follows from Lemma \ref{lemSalzersConstants}
and the bound on $\prod_{k = 0} \wlr{x - x_k}$ presented in \cite{SALZER}
\peProof{Lemma}{lemWarmUp}

\pbProof{Lemma}{lemLp}
Equation \pRef{def_akj} yields
\begin{equation}
\label{zero_sum}
\sum_{j = 0}^n a_{k,j} = \sum_{j = 0}^n z_j \wplagr{j}{c_k} - \wpit{c_k}{z} \sum_{j = 0}^n
\wplagr{j}{c_k} =
\wpit{c_k}{z} - \wpit{c_k}{z} = 0.
\end{equation}
Consider now $\upsilon_j := y_j - y_{j-1}$. It follows that
$y_j = y_0 + \sum_{i = 1}^j \upsilon_i$. Using \pRef{zero_sum}, we obtain
\begin{equation}
\label{hk_small}
\wfc{h_k}{\wvec{y}} = \sum_{j = 0}^n a_{k,j} y_j = a_{k,0} y_0 + \sum_{j = 1}^n a_{k,j} \wlr{y_0 + \sum_{i = 1}^j \upsilon_i}
= \sum_{j = 1}^n a_{k,j} \sum_{i = 1}^j \upsilon_i = \sum_{i = 1}^n \upsilon_i \sum_{j = i}^n a_{k,j}.
\end{equation}
The constraints \pRef{constraint_y} imply that
$\wabs{\upsilon_i} \leq M \wlr{\hat{x}_i - \hat{x}_{i-1}}$.
As a result,  \pRef{hk_small} leads to
\[
\wabs{\wpe{c_k}{y}} =
\wabs{\wfc{h_k}{\wvec{y}}}
\leq \sum_{i = 1}^n \wabs{\upsilon_i} \wabs{\sum_{j = i}^n a_{k,j}} \leq
M \, \sum_{i = 1}^n \wlr{\hat{x}_{i} - \hat{x}_{i-1}} \wabs{\sum_{j = i}^n a_{k,j}}
\]
and we are done
\peProof{Lemma}{lemLp}


\pbProof{Lemma}{lemBoundTau}
We show that for every $\hat{x} \in [\hat{x}^-,\hat{x}^+]$ there exists $x \in [x^-,x^+]$ such that
\begin{equation}
\label{bound_taut}
\wabs{\wplagrxu{k}{\hat{x}}{\hat{x}} \wlr{\hat{x} - \hat{x}_k}} \leq
\frac{\wabs{\wplagrx{k}{x} \wlr{x - x_k}}}{\wlr{1 - \wnorm{\wvec{z}}_\infty}\wlr{1 - \Delta}}.
\end{equation}
When $\hat{x} = \hat{x}_k$ we can satisfy \pRef{bound_taut} by taking $x = x_k$.
For $\hat{x} \not \in \wset{\hat{x}_0,\hat{x}_1,\dots,\hat{x}_n}$,
equations  \pRef{def_zk},  and \pRef{def_lagrange}  and  Lemma \ref{lemBackwardProduct} lead to
\[
\wabs{\frac{\wplagrxu{k}{\hat{x}}{\hat{x}} \wlr{\hat{x} - \hat{x}_k}}{\wplagrx{k}{x} \wlr{x - x_k}}} =
\wabs{\frac{\wfc{\lambda_k}{\wvec{\hat{x}}}}{\wfc{\lambda_k}{\wvec{x}}} \prod_{k = 0}^n \frac{\hat{x} - \hat{x}_k}{x - x_k}}
= \frac{1 + \beta}{1 + z_k},
\]
with $\beta $ such that $\wabs{\beta} \leq \Delta / \wlr{1 - \Delta}$.
Equation \pRef{bound_taut} follows from this equation and we are done
\peProof{Lemma}{lemBoundTau}

\pbProof{Lemma}{lemSterbenz}
Lemma \ref{lemSterbenz} is Theorem 11 in page 38 of \cite{GOLDBERG}
\peProof{Lemma}{lemSterbenz}

\pbProof{Lemma}{lemBoundProd}
If $v_j v_k > 0$ then $\tilde{v}_j = v_j + v_k$ and $\tilde{v}_k = 0$ satisfy
$\wabs{\tilde{v}_j} + \wabs{\tilde{v}_k} = \wabs{v_j + v_k} = \wabs{v_j} + \wabs{v_k}$
and
\[
\wlr{1 + \tilde{v}_j } \wlr{1 + \tilde{v}_k } =
1 + v_j + v_k < 1 + v_j + v_k + v_j v_k = \wlr{1 + v_j} \wlr{1 + v_k}.
\]
Therefore, by replacing $v_j$ and $v_k$ by $\tilde{v}_j$ and $\tilde{v}_k$ we do not change
the sums $\wfc{s_-}{\wvec{v}}$ and $\wfc{s_+}{\wvec{v}}$, and we decrease the
product $\prod_{i = 0}^n \wlr{1 + v_i}$, because the hypothesis $s_- < 1$ implies that
$v_i > -1$ and all its factors are positive.
Applying this argument while there are pairs $v_j$, $v_k$ with $v_j v_k > 0$, we
conclude that
\begin{equation}
\label{magic}
\prod_{k = 0}^n \wlr{1 + v_k} \geq  \wlr{1 - s_-}\wlr{1 + s_+}.
\end{equation}
This equation leads to
\[
\prod_{k = 0}^n \wlr{1 + v_k} - 1 - s \geq - s_- s_+
\]
and the identity $s_- s_+ = (a^2 - s^2)/4 $ yields the first inequality in \pRef{lower_bound_prod}.
Equation \pRef{magic} also shows that
\[
\prod_{k=0}^n \frac{1}{1 + v_k}  - 1 \leq \frac{1}{\wlr{1 - s_-}\wlr{1 + s_+}} - 1 = \frac{s_-}{\wlr{1 - s_-}\wlr{1 + s_+}} - \frac{s_+}{1 + s_+},
\]
and this proves the second inequality in \pRef{bound_inv_prod}.
The set
$\wcal{C}:= \wset{\wvec{v} \in \wrn{n + 1} \ \wrm{with} \ \wfc{s_-}{\wvec{v}} < 1}$ is convex
because the function $s_-$ is convex.
Let us define $h: \wcal{C}\rightarrow \wrone{}$ by
\[
\wfc{h}{\wvec{v}} := \prod_{j = 1}^n \frac{1}{1 + v_j} - \wlr{1 - s}.
\]
The function $h$ has partial derivatives
\[
\wdfdx{h}{v_k}{\wvec{v}} = 1 - \frac{1}{1 + v_k} \prod_{j = 1}^n \frac{1}{1 + v_j},
\]
\[
\wdfdxx{h}{v_k}{\wvec{v}} =  \frac{2}{\wlr{1 + v_k}^2} \prod_{j = 1}^n \frac{1}{1 + v_j}
\hspace{0.7cm} \wrm{and, \ \ for } \ j \neq k, \ \hspace{0.7cm}
\wdfdxy{h}{v_j}{v_k}{\wvec{v}} =  \frac{1}{\wlr{1 + v_j} \wlr{1 + v_k}} \prod_{i = 1}^n \frac{1}{1 + v_i}.
\]
Therefore, its Hessian can be written as
\[
\whessf{h}{\wvec{v}} = \wlr{\prod_{j=1}^n \frac{1}{1 + v_j}} \ \ \wvec{D} \wlr{\wvec{I} + \wvone{} \wvone{}^t} \wvec{D},
\]
where $\wvec{I}$ is the identity matrix, $\wvone{}$ is the vector with all entries equal to $1$
and $\wvec{D}$ is the diagonal matrix with $d_{jj} = \wlr{1 + v_j}^{-1}$. Therefore,
$\whessf{h}{\wvec{v}}$
is positive definite and $h$ is convex. Since $\wfc{h}{0} = 0$
and $\wgradf{h}{0} = 0$ we have that $\wfc{h}{\wvec{v}} \geq 0$ for all $\wvec{v} \in \wcal{C}$.
As a result, we have the first inequality in \pRef{bound_inv_prod}. When $s < 1$, this inequality
leads to \pRef{upper_bound_prod}.

Finally, when $a < 1$ , the inequality $1 - x \leq 1/(1 + x)$ for $x \in (-1,1)$,
and the bound \pRef{magic} lead to
\[
 P \geq  \prod_{k = 1}^n \wlr{1 - \wabs{v_k}} \geq
 \wlr{\, 1 - \wfc{s_-}{\wvec{- \wabs{\wvec{v}}}} \, } \, \wlr{\, 1 + \wfc{s_+}{\wvec{-\wabs{\wvec{v}}}} \, } = 1 - a
 \]
 and
 \[
 P \leq  \prod_{k = 0}^{n} \frac{1}{1 - \wabs{v_k}} \leq
 \frac{1}{\wlr{\, 1 - \wfc{s_-}{- \wvec{\wabs{\wvec{v}}}} \, } \, \wlr{\, 1 + \wfc{s_+}{-\wvec{\wabs{\wvec{v}}}} \, }}
 = \frac{1}{1 - a}.
\]
The bound  $\wabs{P - 1} \leq a / (1 - a)$ follows from the last two equations and we are done
\peProof{Lemma}{lemBoundProd}


\pbProof{Lemma}{lemSums}
Since $x_k := - \cos \frac{k \pi}{n}$, the
identity $\cos a - \cos b = 2 \wfc{\sin}{\frac{a+b}{2}} \wfc{\sin}{\frac{b- a}{2}}$ leads to
\begin{equation}
\label{ratio_rk_a}
\frac{1}{x_k - x_j} = \frac{1}{2 \wfc{\sin}{\frac{k-j}{2n} \pi} \wfc{\sin}{\frac{k+j}{2n} \pi}}.
\end{equation}
If $0 \leq j < k$, then, by the concavity of the $\sin$ function in $[0,k \pi/n] \subset [0,\pi/2]$, we have
\[
\wfc{\sin}{\frac{k-j}{2n} \pi} \geq
\frac{\wfc{\sin}{\frac{k \pi}{2n}}}{\frac{k \pi}{2n}} \times \frac{k - j}{2n} \pi
=
\wfc{\sin}{\frac{k \pi}{2n}} \times \frac{k - j}{k}
\]
and
\[
\wfc{\sin}{\frac{k+j}{2n} \pi} \geq
\frac{\wfc{\sin}{\frac{k \pi}{n} }}{\frac{k \pi}{n}} \times \frac{k + j}{2n} \pi =
\wfc{\sin}{\frac{k \pi}{n} }\times \frac{k + j}{2k}.
\]
Combining the last two inequalities with the well known inequality $\sum_{j = 1}^m \frac{1}{j} \leq \wfc{\log}{2m + 1}$ we obtain
\[
\sum_{j = 0}^{k-1} \frac{1}{x_k - x_j} \leq  \sum_{j = 0}^{k-1}
\frac{k^2}{\wfc{\sin}{\frac{k \pi}{n}} \wfc{\sin}{\frac{k \pi}{2n}} \wlr{k^2 - j^2}}
\]
\[
=
\frac{k}{\wfc{2 \sin}{\frac{k \pi}{n}} \wfc{\sin}{\frac{k \pi}{2n}}} \sum_{j = 0}^{k-1} \wlr{\frac{1}{k - j} + \frac{1}{k+j}} =
\frac{k}{\wfc{2 \sin}{\frac{k \pi}{n}} \wfc{\sin}{\frac{k \pi}{2n}}} \wlr{ \sum_{j = 1}^{k} \frac{1}{j} + \sum_{i=k}^{2k - 1} \frac{1}{j}} =
\]
\[
= \frac{k}{\wfc{2 \sin}{\frac{k \pi}{n}} \wfc{\sin}{\frac{k \pi}{2n}}} \wlr{\frac{1}{k} + \sum_{j = 1}^{2 k - 1} \frac{1}{j}}
\leq \frac{k}{\wfc{2 \sin}{\frac{k \pi}{n}} \wfc{\sin}{\frac{k \pi}{2n}}} \wlr{\frac{1}{k} + \wfc{\log}{4k - 1}}.
\]
This proves the first inequality in \pRef{good_rk_b}. The second inequality in \pRef{good_rk_b} follows
from the first inequality and the observation that the hypothesis $1 \leq k \leq n/2$ and the concavity of $\sin$ in $[0,\pi/2]$ yield
\[
\wfc{\sin}{\frac{k \pi}{n}} \geq \frac{\sin \frac{\pi}{2}}{\frac{\pi}{2}} \frac{k \pi}{n} = \frac{2k}{n}
\hspace{1cm} \wrm{and} \hspace{1cm}
\wfc{\sin}{\frac{k \pi}{2n}} \geq \frac{\sin \frac{\pi}{4}}{\frac{\pi}{4}} \frac{k \pi}{2n} = \frac{\sqrt{2} \, k}{n}.
\]
If $0 \leq k < j \leq n$ then $\wlr{j - k}/\wlr{2n} \pi \leq \pi/2$ and
the concavity of $\sin x$ in $[0, 3 \pi/ 4]$ yields
\begin{equation}
\label{conc_sin}
\wfc{\sin}{\frac{j-k}{2n} \pi} \geq
\frac{\wfc{\sin}{\frac{\pi}{2}}}{\frac{\pi}{2}} \times \frac{j - k}{2n} \pi
=
\frac{j - k}{n}
\hspace{0.2cm} \wrm{and} \hspace{0.2cm}
\wfc{\sin}{\frac{j+k}{2n} \pi} \geq
\frac{\wfc{\sin}{\frac{3 \pi}{4} }}{\frac{3 \pi}{4}} \times \frac{j + k}{2n} \pi =
\sqrt{2} \, \frac{j + k}{3n}.
\end{equation}
If $k = 0$, then equations \pRef{ratio_rk_a} and \pRef{conc_sin} yield
\[
\sum_{j = 1}^n \frac{1}{x_j - x_0} = \sum_{j = 1}^n \frac{1}{2 \wfc{\sin^2}{\frac{j \pi}{2n}}}
\leq \sum_{j = 1}^n \frac{n^2}{2 j^2} < \frac{n^2}{2} \sum_{j=1}^\infty \frac{1}{j^2} = \frac{\pi^2 n^2}{12} \leq 0.82247 n^2.
\]
When $1 \leq k < j \leq n$, equation \pRef{ratio_rk_a} and the inequalities \pRef{conc_sin} lead to
\[
\sum_{j = k+1}^{n} \frac{1}{x_j - x_k} \leq  \frac{3 n^2}{2 \sqrt{2}} \sum_{j = k+1}^{n} \frac{1}{j^2 - k^2} =
\frac{3n^2}{4 \sqrt{2}\, k} \sum_{j = k+1}^{n} \wlr{\frac{1}{j - k} - \frac{1}{j+k}} =
\]
\[
=
\frac{3n^2}{4 \sqrt{2}\, k} \wlr{ \sum_{j = 1}^{n - k} \frac{1}{j} - \sum_{j=2 k + 1}^{n + k} \frac{1}{j}}
\leq
\frac{3n^2}{4 \sqrt{2}\, k} \sum_{j = 1}^{2k} \frac{1}{j}
\leq
\frac{3n^2}{4 \sqrt{2}\, k} \wfc{\log}{4 k + 1}.
\]
Finally, by the concavity of $\wfc{f}{x} = \sin x$ in $[0,\pi/{20}]$, for every $0 \leq i < j \leq n$ we have
\[
\frac{1}{\wabs{x_i - x_j}} \leq \frac{1}{\wabs{x_1 - x_0}} \leq \frac{1}{1 - \wfc{\cos}{\frac{\pi}{n}}} =
\frac{1}{2 \wfc{\sin^2}{\frac{\pi}{2n}}} \leq
\frac{1}{2} \wlr{ \frac{ \sin \frac{\pi}{20}}{\frac{\pi}{20}} \frac{\pi}{2n} }^{-2} \leq 0.20432 n^2,
\]
and this proof is complete
\peProof{Lemma}{lemSums}


\pbProof{Lemma}{lemDelta}
Definition \pRef{def_delta} states that $\delta_{jk} = 0$ when
$j = k$. For $0 \leq j \neq k \leq n$ we have
\[
\wabs{\delta_{jk}} = \wabs{\frac{x_j - x_k}{\hat{x}_j - \hat{x}_k} - 1} =
 \wabs{\frac{\wlr{\hat{x}_j - x_j} - \wlr{\hat{x}_k - x_k}}{\hat{x}_j - \hat{x}_k}} \leq
\frac{2 \wsupdx}{\wabs{x_j - x_k} - 2 \wsupdx}
\]
\[
= \frac{2 \wsupdx}{1 - 2 \wsupdx / \wabs{x_j - x_k}}
\frac{1}{\wabs{x_j - x_k}} \leq \frac{\kappa}{\wabs{x_j - x_k}},
\]
and the first equation in \pRef{bound_delta} holds.
When $j < k$ we have
\[
\max \wset{\wabs{\delta_{jk}},\wabs{\delta_{j\wlr{k+1}}}} \leq
 \frac{\kappa}{x_k - x_j},
\]
when $j > k + 1$,
\[
\max \wset{\wabs{\delta_{jk}},\wabs{\delta_{j\wlr{k+1}}}} \leq
\frac{\kappa}{x_j - x_{k+1}},
\]
and when $j \in \wset{k,k+1}$,
\[
\max \wset{\wabs{\delta_{jk}},\wabs{\delta_{j\wlr{k+1}}}} = \wabs{\delta_{k\wlr{k+1}}} \leq  \frac{\kappa}{x_{k+1} - x_k},
\]
and equation \pRef{bound_big_delta} follows from the last four equations
\peProof{Lemma}{lemDelta}


\pbProof{Lemma}{lemLipschitz}
Let $Q$ be the polynomial
of degree $n$ defined by
\[
\wfc{P}{x;\wvec{x},\wvec{y}} = y_k + \wlr{x - x_k} \wfc{Q}{x}.
\]
For $j \neq k$ we have $\wfc{Q}{x_j} = \wlr{y_j - y_k}/\wlr{x_j - x_k}$ and
Lemma \ref{lemLipschitz} follows from the definition of Lebesgue constant
in \pRef{def_lebesgue}
\peProof{Lemma}{lemLipschitz}

\subsection{Proofs of the theorems}
\pbProof{Theorem}{thmBestApprox}
Note that, for $0 \leq k \leq n$,
\[
\wabs{\wfc{Q}{x_k} - \wfc{P}{x_k;\wvec{x},\wvec{y}}} \leq
 \wabs{\wfc{Q}{x_k} - \wfc{f}{x_k}} + \wabs{\wfc{f}{x_k} - \wfc{P}{x_k;\wvec{x},\wvec{y}}} \leq M,
\]
and equation \pRef{lebesgue_prop} leads to
\begin{equation}
\label{pf_minus_p}
\wabs{\wfc{Q}{x} - \wfc{P}{x;\wvec{x},\wvec{y}}} \leq \Lambda_{\hat{x}^-,\hat{x}^+,\wvec{x}} M
\end{equation}
for all $x \in [\hat{x}^-,\hat{x}^+]$. Similarly,
\[
\wabs{\wfc{Q}{\hat{x}_k} - \wfc{P}{\hat{x}_k;\wvec{\hat{x}},\wvec{y}}} \leq
 \wabs{\wfc{Q}{\hat{x}_k} - \wfc{f}{\hat{x}_k}} + \wabs{\wfc{f}{\hat{x}_k} - \wfc{f}{x_k}}
+
\wabs{\wfc{f}{x_k} - \wfc{P}{\hat{x}_k;\wvec{\hat{x}},\wvec{y}}}
\leq M + L \wabs{\hat{x}_k - x_k},
\]
and
\begin{equation}
\label{pf_minus_p_hat}
\wabs{\wfc{Q}{x} - \wfc{P}{x;\wvec{\hat{x}},\wvec{y}}} \leq
\Lambda_{\hat{x}^-,\hat{x}^+,\wvec{\hat{x}}} \wlr{M + L \wsupdx}.
\end{equation}
Equation \pRef{best_approx} follows from \pRef{pf_minus_p}, \pRef{pf_minus_p_hat} and the triangle inequality.
Equations \pRef{bound_lebesgue_salzer} and  \pRef{mascheroni} lead to
\begin{equation}
\label{salzer_approx_b}
\Lambda_{-1,1,\wvec{\hat{x}}^c} + \Lambda_{-1,1,\wvec{x}^c} \leq 2.0629 \Lambda_{-1,1,\wvec{x}^c} \leq 1.3133 \wlr{\log n + 1.5127}.
\end{equation}
Equation \pRef{bound_salzer_approx} follows from the bound last two bounds and \pRef{best_approx}
\peProof{Theorem}{thmBestApprox}


\pbProof{Theorem}{thmFirstBack}
Given $\hat{x} \in [\hat{x}^-,\hat{x}^+]$,
the arguments used in the proof of Theorem 3.2 in \cite{HIGHAM_STAB} show that
the first formula $p$ with rounded nodes $\wvec{\hat{x}}$  satisfies
\begin{equation}
\label{aux_first_back}
\wrounde{\wfc{p}{\hat{x},\wvec{\hat{x}},\wvec{y},\hat{\wvec{w}}}} =
\sum_{k = 0}^n \hat{w}_k y_k \wst{3 n + 5}_k \prod_{j \neq k} \wlr{\hat{x} - \hat{x}_k},
\end{equation}
where the $\wst{m}$ are the Stewart's relative error counters described in \cite{HIGHAM}.
Lemma \ref{lemBackwardProduct} yields $x$ satisfying \pRef{dx_first_back}
and $\beta_0, \dots, \beta_n$ such that $\wabs{\beta_k} \leq \Delta / \wlr{ 1 - \Delta}$
and
\[
\prod_{j \neq k} \wlr{\hat{x} - \hat{x}_k}  = \wlr{1 + \beta_k} \prod_{j \neq k} \wlr{x - x_k}.
\]
Combining this equation with \pRef{aux_first_back}  we obtain
\begin{equation}
\label{aux_first_back_b}
\wrounde{\wfc{p}{\hat{x},\wvec{\hat{x}},\wvec{y},\hat{\wvec{w}}}} =
\sum_{k = 0}^n w_k y_k \wlr{1 + \beta_k} \wst{3 n + 5}_k \wlr{1 + \frac{\hat{w}_k - w_k}{w_k}} \prod_{j \neq k} \wlr{x - x_k}
= \wfc{p}{x,\wvec{x},\tilde{\wvec{y}},\wvec{w}}
\end{equation}
for $\tilde{y}_k$ in \pRef{y_tilde}
and Theorem \ref{thmFirstBack} follows from Lemma 3.1 in \cite{HIGHAM},
which states that $\wabs{\wst{m} - 1} \leq m \epsilon / \wlr{1 - m \epsilon}$ when $m \epsilon < 1$.
\peProof{Theorem}{thmFirstBack}


\pbProof{Theorem}{thmFirstForward}
Taking $\wvec{x} = \hat{\wvec{x}}$,
and using equations \pRef{first_formula}, \pRef{def_zk} and \pRef{aux_first_back_b},
with $w_k = \wfc{\lambda_k}{\hat{\wvec{x}}}$ and $\beta_k = 0$, and
the fact that the first formula with weights $\wfcc{\lambda}{\wvec{\hat{x}}}$
interpolates $\wvec{y}$ at the nodes $\wvec{\hat{x}}$, we obtain
\[
\wrounde{\wfc{p}{\hat{x},\wvec{\hat{x}},\wvec{y},\hat{\wvec{w}}}} - \wpit{\hat{x}}{y} =
\wlr{\prod_{k = 0}^n \wlr{\hat{x} - \hat{x}_k}} \wlr{
\sum_{k=0}^n \frac{\wfcc{\lambda_k}{\wvec{\hat{x}}} y_k \wlr{1+ z_k} \wst{3n + 5}_k}{\hat{x} - \hat{x}_k} -
\sum_{k=0}^n \frac{\wfcc{\lambda_k}{\wvec{\hat{x}}} y_k}{\hat{x} - \hat{x}_k}
}
\]
\begin{equation}
\label{pre_bound_first}
=
\wlr{\prod_{k = 0}^n \wlr{\hat{x} - \hat{x}_k}} \wlr{
\sum_{k = 0}^n \frac{\wfcc{\lambda_k}{\wvec{\hat{x}}} y_k \nu_k }{\hat{x} - \hat{x}_k}
} = \wpit{\hat{x}}{\wvec{y \bm{\nu}}},
\end{equation}
for
\[
\nu_k := \wlr{1 + z_k} \wst{3 n + 5}_k - 1 = z_k \wst{3 n + 5}_k + \wst{3 n + 5}_k - 1.
\]
The definition of $\wst{3 n + 5}_k$ in \cite{HIGHAM} and Lemma \ref{lemBoundProd}
show that
\[
1 - \wlr{3 n + 5} \epsilon \leq \wst{3 n + 5}_k  \leq \frac{1}{1 - \wlr{3 n + 5} \epsilon},
\]
and lead to \pRef{bound_nuk}
and \pRef{bound_first} follows from \pRef{lebesgue_prop}.
Finally, \pRef{lower_bound_first} follows from
\pRef{pre_bound_first} and Lemma \ref{lemLipschitz}
\peProof{Theorem}{thmFirstForward}


\pbProof{Theorem}{thmMain}
Let us write $x := \hat{x}$, consider the function
\[
\wfc{g}{x,\wvec{y},\wvec{w}} := \sum_{j = 0}^n \frac{w_k y_k}{x - \hat{x}_k},
\]
and $\tilde{\wvec{w}} := \wfc{\lambda_k}{\wvec{\hat{x}}}$. By the definition of $z_k$ in \pRef{def_zk},
$\hat{\wvec{w}} = \tilde{\wvec{w}} + \tilde{\wvec{w}}\wvec{z}$, where $\wvec{\tilde{w}z}$ is the
component-wise product of $\wvec{\tilde{w}}$ and $\wvec{z}$, i.e., $(\tilde{w}z)_k = \tilde{w}_k \, z_k$.
It follows that, for $\wvec{u} \in \wrn{n+1}$,
\[
\wfc{g}{x,\wvec{u},\hat{\wvec{w}}}
= \sum_{j = 0}^n \frac{\hat{w}_k u_k}{x - u_k} =
\sum_{j = 0}^n \tilde{w}_k \frac{\frac{\hat{w}_k - \tilde{w}_k}{\tilde{w}_k} u_k}{x - \hat{x}_k} +
\sum_{j = 0}^n \frac{\tilde{w}_k  u_k}{x - \hat{x}_k} =
\sum_{j = 0}^n \frac{\tilde{w}_k z_k u_k}{x - \hat{x}_k} + \sum_{j = 0}^n \frac{\tilde{w}_k  u_k}{x - \hat{x}_k}
= \wfc{g}{x,\wvec{uz},\tilde{\wvec{w}}} +\wfc{g}{x,\wvec{u},\tilde{\wvec{w}}},
\]
and we can write the difference
$\wfc{q}{x,\wvec{\hat{x}},\wvec{y},\hat{\wvec{w}}} - \wfc{q}{x,\wvec{\hat{x}},\wvec{y},\tilde{\wvec{w}}}$
in \pRef{good_news} as
\[
\wfc{S}{x,\wvec{y},\wvec{z}} :=
\frac{\wfc{g}{x,\wvec{y},\hat{\wvec{w}}}}
     {\wfc{g}{x,\wvone{},\hat{\wvec{w}}}} -
\frac{\wfc{g}{x,\wvec{y},\tilde{\wvec{w}}}}
     {\wfc{g}{x,\wvone{},\tilde{\wvec{w}}}}
=
\frac{\wfc{g}{x,\wvec{y},\tilde{\wvec{w}}} + \wfc{g}{x,\wvec{yz},\tilde{\wvec{w}}}}
     {\wfc{g}{x,\wvone{},\tilde{\wvec{w}}} + \wfc{g}{x,\wvec{z},\tilde{\wvec{w}}}} -
\frac{\wfc{g}{x,\wvec{y},\tilde{\wvec{w}}}}
     {\wfc{g}{x,\wvone{},\tilde{\wvec{w}}}} = \frac{N + \Delta \! N}{D + \Delta \! D} - \frac{N}{D}
\]
for
$N := \wfc{g}{x,\wvec{y},\tilde{\wvec{w}}}$, $\Delta \! N := \wfc{g}{x,\wvec{yz},\tilde{\wvec{w}}}$,
$D := \wfc{g}{x,\wvone{},\tilde{\wvec{w}}}$, and $\Delta \! D := \wfc{g}{x,\wvec{z},\tilde{\wvec{w}}}$.
Since $\frac{\Delta \! N}{D} - \frac{N}{D} \frac{\Delta \! D}{D}$ is equal to the Error Polynomial $
\wpe{x}{y}$ in \pRef{def_error_pol}, we get
\begin{equation}
\label{aux_main}
\wfc{S}{x,\wvec{y},\wvec{z}} = \frac{1}{D}
\wlr{ \frac{N + \Delta \! N}{1 + \frac{\Delta \! D}{D}} - N }
=  \frac{1}{{1 + \frac{\Delta \! D}{D}}} \wpe{x}{y}.
\end{equation}
The ratio $\frac{\Delta \! D}{D}$ is equal to $\wfc{P}{x, \hat{\wvec{x}}, \wvec{z}}$.
Therefore, the denominator of \pRef{aux_main} satisfies
\begin{equation}
\label{main_at_last}
1 + \Lambda_{\hat{x}^-,\hat{x}^+,\wvec{\hat{x}}} \wnorm{\wvec{z}}_\infty \geq
1 + \wabs{\frac{\Delta \! D}{D}} \geq
\wabs{1 + \frac{\Delta \! D}{D}} \geq 1 - \wabs{\frac{\Delta \! D}{D}} = 1 - \wabs{\wfc{P}{x,\hat{\wvec{x}},\wvec{z}}}
\geq 1 - \Lambda_{\hat{x}^-,\hat{x}^+,\wvec{\hat{x}}}\wnorm{\wvec{z}}_\infty.
\end{equation}
The forward bound \pRef{good_news} follows from \pRef{aux_main} and \pRef{main_at_last}
 and we are done.
\peProof{Theorem}{thmMain}


\pbProof{Theorem}{thmLipschitz}
Expanding the polynomials $\wpit{x}{yz}$ in
\pRef{def_error_pol} in Lagrange's basis we obtain
\[
\wabs{\wpe{x}{y}} =
\wabs{\sum_{k=0}^n z_k \wplagr{k}{x} \wlr{y_k - \wpit{x}{y} }} \leq
\]
\[
\leq \sum_{k=0}^n \wabs{z_k \wplagr{k}{x}}
\wabs{\wfc{f}{\hat{x}_k} - \wfc{f}{x}} +
\sum_{k=0}^n \wabs{z_k \wplagr{k}{x}}
 \wabs{\wfc{f}{x} - \wpit{x}{y}} \leq
 \]
 \[
L \wlr{\max_{0 \leq k \leq n} \wabs{ \wplagr{k}{x} \wlr{x - \hat{x}_k}}} \,
 \sum_{k=0}^n \wabs{z_k}  +
\Lambda_{\hat{x}^-,\hat{x}^+,\wvec{\hat{x}}}
 \wnorm{\wvec{z}}_\infty  \sup_{\hat{x}^- \leq x \leq \hat{x}^+} \wabs{\wfc{f}{x} - \wfc{P}{x,\wvec{\hat{x}},\wvec{y}}}.
\]
This inequality yields \pRef{lipschitz_bound_a}.
In order to prove \pRef{lipschitz_bound_b}, note that, according to
Theorem 16.5 in page 196 of \cite{POWELL} (Jackson's Theorem),
there exists a polynomial $P^*$ of degree $n$ such that
\begin{equation}
\label{powell}
\sup_{x \in [\hat{x}^-,\hat{x}^+]} \wabs{\wfc{f}{x} - \wfcc{P^*}{x}} \leq \frac{L \wlr{\hat{x}^+ - \hat{x}^-} \pi}{4 \wlr{n + 1}}.
\end{equation}
Using the well known bound
\[
\max_{\hat{x}^- \leq x \leq \hat{x}^+} \wabs{\wfc{f}{x} - \wfc{P}{x,\wvec{\hat{x}},\wvec{y}}}_\infty \leq
\wlr{1 + \Lambda_{\hat{x}^-,\hat{x}^+,\wvec{\hat{x}}}} \sup_{\hat{x}^- \leq x \leq \hat{x}^+} \wabs{\wfc{f}{x} - \wfcc{P^*}{x}}
\]
and \pRef{powell} we deduce the that
\[
 \Lambda_{\hat{x}^-,\hat{x}^+,\wvec{\hat{x}}} \wnorm{\wvec{z}}_\infty  \sup_{\hat{x}^- \leq x \leq \hat{x}^+} \wabs{\wfc{f}{x} - \wfc{P}{x,\wvec{\hat{x}},\wvec{y}}}
 \leq
\Lambda_{\hat{x}^-,\hat{x}^+,\wvec{\hat{x}}} \wlr{1 + \Lambda_{\hat{x}^-,\hat{x}^+,\wvec{\hat{x}}}} \wnorm{\wvec{z}}_\infty
\frac{L \wlr{\hat{x}^+ - \hat{x}^-} \pi}{4 \wlr{n + 1}},
\]
and \pRef{lipschitz_bound_b} follows from this last bound and \pRef{lipschitz_bound_a}
\peProof{Theorem}{thmLipschitz}

\pbProof{Theorem}{thmResumeLip}
The definition of Lagrange polynomial \pRef{def_lagrange} leads to
\[
\wplagrxu{k}{x}{x} \wlr{x - x_k} = \wfc{\lambda_k}{\wvec{x}} \prod_{k = 0}^n \wlr{x - x_k}.
\]
In \cite{SALZER}'s notation, we have $\wfc{\lambda_k}{\wvec{x}^c} = 1/\wdfc{\phi_{n+1}}{x_k^c}$ and from its equations (5)
and (6) we obtain $\wfc{\lambda_k}{\wvec{x}^c} \leq  2^{n-1}/n$. At the top of the second column in page 156 of \cite{SALZER} we
learn that $\wabs{\prod_{k = 0}^n \wlr{x - x^c_k}} \leq 2^{1 - n}$ for $x \in [-1,1]$, and combining this bounds we obtain that
\begin{equation}
\label{salzer_tau}
\wfc{\tau}{\wvec{x}^c} = \max_{x \in [-1,1]} \wabs{\wplagrxu{k}{x}{x^c} \wlr{x - x^c_k}} \leq \wlr{2^{n-1}/n} \times 2^{1 - n} =  1/n.
\end{equation}
The bounds  \pRef{sharpZ}, \pRef{sharp_big_delta}, \pRef{sharpZ1} and \pRef{salzer_tau} and Lemma \ref{lemBoundTau}
lead to
\[
\wfc{\tau}{\wvec{\hat{x}}^c} \wnorm{\wvec{z}}_1 \leq \frac{1.0097}{n} \times 3.2704 \times 10^{-16} n^2 \wlr{2.9 + \log n}^2
\leq 3.3022 \times 10^{-16} n \wlr{2.9 + \log n}^2.
\]
The bounds \pRef{sharpZ} and \pRef{bound_lebesgue_salzer} show that, in Salzer's case,
for $A = 0.67667$ and $B = 1.0236$,
\[
\Lambda_{\hat{x}^-,\hat{x}^+,\wvec{\hat{x}}} \wlr{1 + \Lambda_{\hat{x}^-,\hat{x}^+,\wvec{\hat{x}}}}\wnorm{\wvec{z}}_\infty \frac{\wlr{\hat{x}^+ - \hat{x}^-} \pi}{4 \wlr{n + 1}}
\leq  \wlr{A \log n + B} \wlr{A \log n + B + 1}  \times 1.1328 \times 10^{-15} n^2 \frac{\pi}{2 \wlr{n + 1}}
\]
\[
\leq \wlr{\log n + B/A} \wlr{\log n + \wlr{B + 1}/A} \times A^2 \times 1.7794 \times 10^{-15} n
\leq 8.1476 \times 10^{-16} n \wlr{2.9 + \log n}^2,
\]
because $\wlr{x + \frac{B}{A}} \wlr{x + \frac{B + 1}{A}} \leq \wlr{2.9 + x}^2$ for all $x \geq 0$.
The last two bounds using Theorem \ref{thmLipschitz} yield
\[
\wabs{\wpeu{\hat{x}}{y}} \leq 1.1450 \times 10^{-15} L n \wlr{2.9 + \log n}^2,
\]
and Lemma \ref{lemSalzersConstants} and Theorem \ref{thmMain} lead to
\[
\wabs{\wfc{q}{\hat{x};\wvec{\hat{x}}^c,\wvec{y}, \wfc{\lambda}{\wvec{x}^c}} - \wfc{P}{\hat{x};\wvec{\hat{x}}^c,\wvec{y}}}
\leq 1.2042 \times 10^{-15} L n \wlr{2.9 + \log n}^2,
\]
and we proved \pRef{bound_resume_lip}. Equation \pRef{bound_resume_lip_b} follows from this bound, $\wnorm{\hat{\wvec{x}} - \wvec{x}}_\infty \leq 4.6\times 10^{-6}$,
$n \geq 10$ and \pRef{bound_salzer_approx}
\peProof{Theorem}{thmResumeLip}

\pbProof{Theorem}{thmEmpirical}
The hypothesis of this theorem asks for $n \leq 10^{6}$.
Therefore, $2 n \leq 2 \times 10^{16}$ and the bounds \pRef{lip_bound} and
\pRef{bound_lebesgue_salzer} for $2 n + 1$ nodes yield
\begin{equation}
\max_{\hat{x}^- \leq x \leq \hat{x}^+} \wabs{\wpe{x}{y}} \leq
b_n L \wlr{A \log n + B + A \log 2}
\end{equation}
for $A = 0.67667$ and $B = 1.0236$ . Theorem \ref{thmMain}, \pRef{sharpZ} and
\pRef{bound_lebesgue_salzer} for $n$ nodes lead to an upper bound of
\begin{equation}
\label{empirical_aux}
\frac{A L b_n }{1 - 10.841 \times 4.5312 \times 10^{-3}} \wlr{\log n + \frac{B}{A} + \log 2}
\leq 0.71163 L b_n \wlr{\log n + 2.2059},
\end{equation}
and this leads to \pRef{bound_empirical}. Finally, \pRef{bound_empirical_b} follows
from \pRef{empirical_aux}, $\wnorm{\hat{\wvec{x}} - \wvec{x}}_\infty \leq 4.6\times 10^{-6}$
and Theorem \ref{thmBestApprox}.
\peProof{Theorem}{thmEmpirical}


\pbProof{Theorem}{thmStepIIISec}
The hypothesis $\wnorm{\wfc{\bm{\zeta}}{\wfc{\lambda}{\hat{\wvec{x}}},\hat{\wvec{w}}}}_\infty \wlr{1 + \Lambda_{\hat{x}^-,\hat{x}^+,\hat{\wvec{x}}}} < 1$
allows us to use theorem 3 in \cite{MascCam} with $d = 0$,
$\wvec{x} = \hat{\wvec{x}}$ and $\wvec{w} = \wfc{\lambda}{\hat{\wvec{x}}}$  and conclude that
the Lebesgue constant $\Lambda_{\hat{x}^-,\hat{x}^+,\hat{\wvec{x}}, \hat{\wvec{w}}}$, which is defined
in the more general context of rational interpolation in \cite{MascCam}, satisfies $\Lambda_{\hat{x}^-,\hat{x}^+,\hat{\wvec{x}}, \hat{\wvec{w}}} \leq \Lambda$.
The hypothesis $\wlr{n + 2} \wlr{2 + \Lambda} \epsilon < 1$ and Theorem 1 in \cite{MascCam} with $\wvec{x} = \hat{\wvec{x}}$ and $\wvec{w} = \hat{\wvec{w}}$
lead to
\[
\wrounde{\wfc{q}{\hat{x}, \hat{\wvec{x}}, \wvec{y}, \hat{\wvec{w}}}} =
\wfc{q}{\hat{x}, \hat{\wvec{x}}, \tilde{\wvec{y}}, \hat{\wvec{w}}} \hspace{1cm} \wrm{with} \hspace{1cm}
\tilde{y}_k = y_k \wlr{1 + \alpha_k} \wlr{ 1 + \nu_k},
\]
with
\[
\wnorm{\bm{\alpha}}_\infty \leq \frac{\wlr{1 + \Lambda} \wlr{n + 2} \epsilon}{1 - \wlr{n + 2} \wlr{2 + \Lambda} \epsilon } = \theta \epsilon
\hspace{1.0cm} \wrm{and} \hspace{1.0cm}
\wnorm{\bm{\nu}}_\infty \leq \frac{ \wlr{2 n + 6} \epsilon}{1 - \wlr{2 n + 6} \epsilon}
\]
It follows that
\[
\wnorm{ \tilde{\wvec{y}} - \wvec{y}}_\infty \leq \frac{\theta + \wlr{2 n + 6} }{1 - \wlr{2 n + 6} \epsilon} \wnorm{\wvec{y}}_\infty \epsilon,
\]
The definition of $\Lambda_{\hat{x}^-,\hat{x}^+,\hat{\wvec{x}}, \hat{\wvec{w}}}$ in \cite{MascCam} leads to
\[
\wabs{\wrounde{\wfc{q}{\hat{x}, \hat{\wvec{x}}, \wvec{y}, \hat{\wvec{w}}
}} - \wfc{q}{\hat{x}, \hat{\wvec{x}}, \tilde{\wvec{y}}, \hat{\wvec{w}}}} \leq \Lambda \wnorm{\tilde{\wvec{y}} - \wvec{y}}_\infty,
\]
and \pRef{bound_stepIII_a} follows from the last two inequalities. In Salzer's case,
Lemmas \ref{lemSalzersConstants} and \ref{lemLebesguePerturbation} show that
\[
\Lambda \leq \frac{0.67667 \log n + 1.0236}{1 - 4.5312 \times 10^{-3} \times 11.841} \leq 0.71504 \times \wlr{1.5128 + \log n}
\leq 11.456,
\]
and noting that $n + 2 \leq (2 n + 6) / 2$ and using this bound on $\Lambda$ we deduce that
\[
\theta \leq 0.35753 \times \wlr{2.9114 + \log n} \wlr{2 n + 6}
\hspace{0.3cm}
\wrm{and}
\hspace{0.3cm}
\wnorm{ \tilde{\wvec{y}} - \wvec{y}}_\infty \leq 0.35754 \times \wlr{5.7084 + \log n} \wlr{2 n + 6}
\wnorm{\wvec{y}}_\infty \epsilon,
\]
and \pRef{bound_stepIII_b} follows from the last two inequalities
\peProof{Theorem}{thmStepIIISec}


\pbProof{Theorem}{thmGoodRik}
We start by showing that for
each $0 \leq k \leq n$ there exists $\epsilon_k$ such that
\begin{equation}
\label{rik_round}
\tilde{x}_k + 1 = \wlr{x_k + 1} \wlr{1 + \epsilon_k} \hspace{1cm} \wrm{and} \hspace{1cm} \wabs{\epsilon_k} \leq \wnorm{\bm{\theta}}_\infty.
\end{equation}
We analyze the $x_k$ in the three bins:
\begin{itemize}
\item If $x_k$ is in the first bin then $\epsilon_k = \theta_k$ is appropriate, because
$\tilde{x}_k + 1 = \wlr{1 + \theta_k} r_k = \wlr{x_k + 1} \wlr{1 + \theta_k}$.
\item If $x_k$ is in the second bin then $2 \wabs{x_k} \leq 1$ and
\[
\tilde{x}_k + 1 = 1 + x_k \wlr{1 + \theta_k} = \wlr{x_k + 1}\wlr{1 + \epsilon_k}
\hspace{1cm} \wrm{for} \hspace{1cm}
\epsilon_k := \frac{x_k \theta_k}{1 + x_k}.
\]
This $\epsilon_k$ satisfies \pRef{rik_round} because $\wabs{x_k /(1 + x_k)} \leq 1$ when $2 \wabs{x_k} \leq 1$.
\item If $x_k$ is in the third bin then $\tilde{x}_k = 1 + \hat{r}_k = 1 + r_k \wlr{1 + \theta_k}$, with
$r_k = x_k - 1$. Thus,
\[
\tilde{x}_k + 1 = 2 + (x_k - 1) \wlr{1 + \theta_k} = \wlr{x_k + 1} \wlr{1 + \epsilon_k}
\hspace{0.2cm} \wrm{for} \hspace{0.2cm} \epsilon_k := \theta_k \frac{x_k - 1}{x_k + 1}.
\]
This $\epsilon_k$ is valid because $0 \leq x_k \leq 1$ implies that
$-1 \leq \frac{x_k - 1}{x_k + 1} \leq 0$.
\end{itemize}
By symmetry, we need to verify \pRef{good_rik} only for $0 \leq k \leq n/2$.
Let us then assume from now on that $0 \leq  k \leq n/2$ and $0 \leq j \neq k \leq n$.
Defining
\[
\delta'_{jk} := \frac{\tilde{x}_j - \tilde{x}_k}{x_j - x_k} - 1,
\]
we obtain
\[
\wabs{\delta'_{jk}} = \frac{ \wabs{\wset{\wlr{\tilde{x}_k + 1} - \wlr{x_k + 1}}  -
\wset{\wlr{\tilde{x}_j + 1} - \wlr{x_j + 1}}}}{\wabs{x_k - x_j }}
=   \wabs{\frac{ \wlr{x_k + 1} \epsilon_k - \wlr{x_j + 1} \epsilon_j}{x_k - x_j}}
\]
\[
= \wabs{\epsilon_j + \frac{ \wlr{\epsilon_k - \epsilon_j} \wlr{x_k + 1}}{x_k  - x_j}}
\leq
\wnorm{\bm{\theta}}_\infty \wlr{1 + 2 \frac{x_k + 1}{\wabs{x_k - x_j}}}.
\]
Since we are assuming that $-1 \leq x_k \leq 0$, it follows from \pRef{boundMu} and \pRef{def_delta} that
\begin{equation}
\label{good_rk_a}
\wabs{\delta_{jk}} = \wabs{\frac{\delta'_{jk}}{1 + \delta'_{jk}}} \leq \psi \wlr{1 + 2 \frac{x_k + 1}{\wabs{x_k - x_j}}}
\end{equation}
for
\begin{equation}
\label{def_psi}
\psi := \frac{\wnorm{\bm{\theta}}_\infty}{1 - \max_{0 \leq j \neq k \leq n} \wabs{\delta'_{jk}}} \leq
\frac{\wnorm{\bm{\theta}}_\infty}{1 - \wnorm{\bm{\theta}}_\infty \wlr{1 + 2 \times 0.20432 \times n^2  }}
\leq 1.0008 \wnorm{\bm{\theta}}_\infty.
\end{equation}
In particular, if $k = 0$ then $x_k = -1$ and this equation shows that $\wabs{\delta_{j0}} \leq \psi$ for all $j$. This
implies that $\sum_{j \neq 0} \wabs{\delta_{j 0}} \leq n \psi$
and \pRef{good_rik} is satisfied for $k = 0$. Therefore, from now on we assume that $k > 0$.

Combining equation \pRef{good_rk_a} with Lemma \ref{lemSums} we obtain
\begin{equation}
\label{almost_there}
\frac{1}{\psi} \sum_{j \neq k} \wabs{\delta_{j k}}
\leq n + 2 \wlr{x_k + 1} \wlr{
\frac{1 + k \, \wfc{\log}{4k - 1}}{2 \wfc{\sin}{\frac{k \pi}{n}} \wfc{\sin}{\frac{k \pi}{2n}} }  +
\frac{3 n^2 }{4 \sqrt{2} \, k} \wfc{\log}{4 k + 1}}.
\end{equation}
Let us now analyze the case $1 \leq k \leq \frac{n}{3}$. The definition $x_k := - \cos \frac{k \pi}{n}$ yields
\[
x_k + 1 =  2 \, \wfc{\sin^2}{\frac{k \pi}{2n}}.
\]
Equation \pRef{almost_there}, the assumption $k \leq n/3$
and the inequality $\wfc{\sin}{\frac{k \pi}{n}} \geq \wfc{\sin}{\frac{k \pi}{2n}}$ yield
\[
\frac{1}{\psi} \sum_{j \neq k} \wabs{\delta_{j k}}
\leq n +
2 \wlr{1 + k \, \wfc{\log}{4k - 1}} + \frac{3 \pi^2 k }{4 \sqrt{2}} \wfc{\log}{4 k + 1}
\]
\[
\leq n + 2 + n \wlr{\frac{2}{3} + \frac{\pi^2}{4 \sqrt{2}} } \wfc{\log}{\frac{4n}{3} + 1}
\leq n + 2 + 2.5 \ n \ \wfc{\log}{\frac{4n}{3} + \frac{4}{3}}
\leq 2 + 1.8 n + 2.5 \ n \ \wfc{\log}{n + 1}
\]
and \pRef{good_rik} holds for $1 \leq k \leq n/3$.
Finally, if $n/3 < k \leq n/2$, then equation \pRef{good_rk_a} and
Lemma \ref{lemSums} lead to
\[
\frac{1}{\psi} \sum_{j \neq k} \wabs{\delta_{j k}}
\leq n + \frac{n^2}{2 \sqrt{2}} \wlr{\frac{1}{k^2} + \frac{\wfc{\log}{4k - 1}}{k}
+ 3 \frac{\wfc{\log}{4k + 1}}{k}
}.
\]
The expression above decreases as $k$ increases for $k \geq n/3$. Therefore,
replacing $k$ by $n/3$ we get
\[
\frac{1}{\psi} \sum_{j \neq k} \wabs{\delta_{j k}}
\leq n +
\frac{n^2}{2 \sqrt{2}} \wlr{\frac{9}{n^2} + 3 \frac{\wfc{\log}{\frac{4n}{3} + 1}}{n}
+ 9 \frac{\wfc{\log}{\frac{4n}{3} + 1}}{n}
}
\]
\[
= n + \frac{9}{2 \sqrt{2}} +
\frac{6 n}{\sqrt{2}} \ \wfc{\log}{\frac{4n}{3} + 1} \leq 3.1820 + n + 4.2427 \, n \, \wfc{\log}{\frac{4n}{3} + \frac{4}{3}}
\]
\[
\leq 3.1820 + 2.2213 n + 4.2427 n \, \wfc{\log}{n+1}.
\]
This equation and  \pRef{def_psi} lead to \pRef{good_rik} and we are done.
\peProof{Theorem}{thmGoodRik}

\section{Experimental details}
\label{section_experiments}
The data in our tables and plots were generated with \verb C++  code compiled
in Ubuntu 12.04 with gcc 4.8.1 with options \verb -m64  \verb -std=c++11 \verb -O3 \verb -mavx \verb -Wall ,
with \verb NDEBUG  defined. The code was executed on a Intel Core i7 2700K processor.
We used the IEEE-754 double precision arithmetic, with \verb C++'s  type \verb double .
These numbers are used by \verb Matlab  and correspond to \verb real*8  in \verb Fortran .

We used gcc 4.8.1's quadruple precision type \verb __float128  as a benchmark:
we considered results obtained using this
arithmetic as exact. We checked our results with Intel's \verb C++  13.0 compiler with
option \verb -Qoption,cpp,--extended_float_type  and
its quadruple precision type \verb _Quad .
We also performed accuracy experiments with the versions of the compilers above for
Windows 7 and OS X, in a Quad Core Intel Xeon processor.
There were no relevant differences in the results.

The sets $X_{-1,n}$ and $X_{0,n}$ in tables \ref{table_rho},
\ref{table_stepII_error} and \ref{table_overall_error} contain $10^5$ points each. These points
are distributed in 100 intervals $\wlr{x_k,x_{k+1}}$. $X_{-1,n}$ uses
$0 \leq k < 100$ and $X_{0,n}$ considers $n/2 - 100 \leq k < n/2$.
In each interval $\wlr{x_k,x_{k+1}}$ we picked the 200 floating point numbers to the right
of $x_k$ and the 200 floating point numbers to the left of $x_{k+1}$. The remaining
600 points where equally spaced in $\wlr{x_k,x_{k+1}}$.
The Step II errors in tables \ref{table_rho} and \ref{table_overall_error}
were estimated by performing Step II in double precision,
Step III was evaluated with gcc 4.8.1's \verb __float128  arithmetic,
 and the result was compared with the interpolated function evaluated with \verb __float128  arithmetic.

The trial points in Figure \ref{figure_least_squares} were chosen as follows: for each $n$
we considered the relative errors $z_k^s$ and $z_k^r$ for the Salzer's weights
and the numerical weights. We then we picked 4 groups of ten indexes:
the indexes corresponding largest values $y_k z_k^s$, the indexes
corresponding to the ten smallest values of $y_k z_k^s$ and the analogous 20 indexes
for $z_k^r$. We then removed the repeated indexes and obtained a vector with
$n_i$ indexes. For each index $k > 0$ we picked
the $2000$ floating points to the left of $\hat{x}_k$ and for each index $k < n$ we
picked the $2000$ floating points to the right of $\hat{x}_k$. We then
considered the $n_i$ intervals of the form $[x_{k-1},x_{k}]$ for $k = 1, n/n_i, 2 n / n_i \dots$
and picked 2000 equally spaced points in each of these intervals.
The first formula was evaluated at these points as
in subsection 3.1 of \cite{Masc}.

The $b_n$ in Table \ref{table_bn} were computed using \verb __float128  arithmetic,
from nodes $\hat{x}_k$ obtained from the formula
$x_k = \wfc{\sin}{\frac{2k - n}{2n} \pi}$ using IEEE-754 double arithmetic.
The versions with $3$ bins in tables \ref{table_stepII_error}, \ref{table_overall_error}
and \ref{table_times} are as in Figure \ref{figure_three_bins}.
The versions with $39$ bins consider the central bin $[-2^{-10},2^{-10}]$, with base $b_{20} = 0$,
the bins $[-1,2^{-10} - 1)$, $[2^{-k} - 1, 2^{1 - k} - 1)$ for $2 \leq k \leq 10$
and $[- 2^{-k}, -2^{-k - 1})$ for $1 \leq k < 9$, with base at their left extreme
point. The remaining 19 bins and bases were obtained by reflection around $0$.
The versions with $79$ bins consider the central bin $[-2^{-20},2^{-20}]$,
with base $b_{40} = 0$, the bins
$[-1,2^{-20} - 1)$, $[2^{-k} -1, 2^{1 - k}-1)$ for $1 \leq k < 20$
and $[-2^{-k}, -2^{-k-1})$ for  $1 \leq k < 19$, with base at their
left extreme point. The remaining 39 bins and bases were obtained by
reflection.

The times in Table \ref{table_times} were measured using the cpu timer available
in version 1.54.0 of the boost library \cite{BOOST_SITE}. This timer
measures the time taken only by the process one is concerned with.

\end{document}